\documentclass{amsart}

\usepackage{amscd,amsmath,xypic,amssymb,combelow,tikz-cd,etoolbox,calligra,mathrsfs,enumitem,mathtools,hyperref}

\makeatletter
\patchcmd{\@settitle}{\uppercasenonmath\@title}{}{}{}
\makeatother

\xyoption{all}
\CompileMatrices

\makeatletter
\@addtoreset{equation}{section}
\makeatother

\newtheorem{theorem}[subsection]{Theorem}

\newtheorem{proposition}[subsection]{Proposition}
\newtheorem{lemma}[subsection]{Lemma}

\newtheorem{definition}[subsection]{Definition}

\newtheorem{claim}[subsection]{Claim}
\newtheorem{example}[subsection]{Example}
\newtheorem{remark}[subsection]{Remark}

\def\loccitt{\emph{loc. cit.}}

\def\fg{{\mathfrak{g}}}

\def\fm{{\mathfrak{m}}}

\def\BA{{\mathbb{A}}}
\def\BC{{\mathbb{C}}}
\def\BK{{\mathbb{K}}}
\def\BL{{\mathbb{L}}}

\def\BN{{\mathbb{N}}}
\def\BF{{\mathbb{F}}}

\def\BQ{{\mathbb{Q}}}
\def\BZ{{\mathbb{Z}}}

\def\CI{{\mathcal{I}}}
\def\CJ{{\mathcal{J}}}

\def\CM{{\mathcal{M}}}
\def\CN{{\mathcal{N}}}

\def\CS{{\mathcal{S}}}

\def\CV{{\mathcal{V}}}
\def\CW{{\mathcal{W}}}

\def\ph{\varphi}

\def\vs{\varsigma}

\def\and{\textrm{ }\&\textrm{ }}

\def\sym{\textrm{sym}}
\def\Sym{\textrm{Sym}}

\def\esym{\emph{sym}}
\def\eSym{\emph{Sym}}

\def\oCS{\mathring{\CS}}

\def\nn{{{\BN}}^I}
\def\zz{{{\BZ}}^I}

\def\UU{\mathbf{U}}
\def\UUo{\mathbf{U}^0}
\def\UUp{\mathbf{U}^+}
\def\UUpm{\mathbf{U}^\pm}
\def\UUm{\mathbf{U}^-}
\def\UUg{\mathbf{U}^\geq}
\def\UUl{\mathbf{U}^\leq}

\def\tUU{{\widetilde{\mathbf{U}}}}
\def\tUUp{{\widetilde{\mathbf{U}}^+}}
\def\tUUpm{{\widetilde{\mathbf{U}}^\pm}}
\def\tUUo{{\widetilde{\mathbf{U}}^0}}
\def\tUUm{{\widetilde{\mathbf{U}}^-}}

\def\bs{{\boldsymbol{\vs}}}

\def\b0{{\boldsymbol{0}}}

\def\ov{{\overline{v}}}
\def\ow{{\overline{w}}}

\def\bk{\boldsymbol{k}}
\def\bn{\boldsymbol{n}}

\def\oij{\overrightarrow{ij}}
\def\oji{\overrightarrow{ji}}

\def\tUpsilon{\widetilde{\Upsilon}}
\def\tzeta{\widetilde{\zeta}}

\def\oBK{\overline{\BK}}

\begin{document}

\title[Quantum loop groups for arbitrary quivers]{\Large{\textbf{Quantum loop groups for arbitrary quivers}}}

\author[Andrei Negu\cb t]{Andrei Negu\cb t}
\address{MIT, Department of Mathematics, Cambridge, MA, USA}
\address{Simion Stoilow Institute of Mathematics, Bucharest, Romania}
\email{andrei.negut@gmail.com}

\maketitle

\begin{abstract} 

We study the dual constructions of quantum loop groups and Feigin-Odesskii type shuffle algebras for an arbitrary quiver, for which the arrow parameters are arbitrary non-zero elements of any field. Examples of our setup include $K$-theoretic Hall algebras of quivers with 0 potential, quantum loop groups of Kac-Moody type and quiver quantum toroidal algebras.

\end{abstract}

\section{Introduction}
\label{sec:intro}

\medskip

\subsection{}
\label{sub:general intro}

In the present paper, we develop a systematic treatment of quantum loop groups, which generalizes the particular cases treated in \cite{QQTA, Loop, Quiver 1, Quiver 3}. Specifically, we fix a finite set $I$, a field $\BK$ of characteristic 0, and a collection of rational functions
\begin{equation}
\label{eqn:general zeta intro}
\zeta_{ij}(x) \in \frac {\BK[x^{\pm 1}]}{(1-x)^{\delta_{ij}}}
\end{equation}
for all $i,j \in I$. To this datum, one can associate two objects

\medskip

\begin{itemize}

\item The (half) \textbf{quadratic quantum loop groups} 
\begin{align}
&\tUUp = \BK \Big\langle e_{i,d} \Big \rangle_{i \in I, d \in \BZ} \Big / \text{relation \eqref{eqn:rel quad}} \label{eqn:quad intro plus} \\
&\tUUm = \BK \Big\langle f_{i,d} \Big \rangle_{i \in I, d \in \BZ} \Big / \text{relation \eqref{eqn:opposite}} \label{eqn:quad intro minus}
\end{align}

\medskip

\item The \textbf{big shuffle algebras} (following ideas of \cite{E,FO})
\begin{equation}
\label{eqn:big shuf intro}
\CV^\pm = \bigoplus_{\bn = (n_i)_{i \in I} \in \nn} \BK [z_{i1}, z_{i1}^{- 1}, \dots, z_{in_i}, z_{in_i}^{- 1}]_{i \in I}^{\sym}
\end{equation}
which are endowed with the product \eqref{eqn:shuf prod} and its opposite, respectively. \medskip

\end{itemize}

\medskip

\noindent The constructions above are related by $\BK$-algebra homomorphisms
\begin{equation}
\label{eqn:upsilon intro}
\tUUpm \xrightarrow{\tUpsilon^\pm} \CV^\pm, \qquad e_{i,d}, f_{i,d} \mapsto z_{i1}^d
\end{equation}
and pairings 
\begin{align}
&\tUUp \otimes \CV^- \xrightarrow{\langle \cdot, \cdot \rangle} \BK \label{eqn:pair intro 1} \\
&\CV^+ \otimes \tUUm \xrightarrow{\langle \cdot, \cdot \rangle} \BK \label{eqn:pair intro 2}
\end{align}
whose formulas will be recalled in \eqref{eqn:pair formula} and \eqref{eqn:pair formula opposite}. With this in mind, define

\medskip

\begin{itemize}

\item The (half) \textbf{quantum loop groups} as
\begin{equation}
\label{eqn:first def}
\UUpm = \tUUpm \Big / K^\pm 
\end{equation}
where $K^\pm = \text{Ker }\tUpsilon^\pm$.

\medskip

\item The \textbf{shuffle algebras} as
\begin{equation}
\label{eqn:second def}
\oCS^\pm = \text{Im } \tUpsilon^\pm \subset \CV^\pm
\end{equation}

\end{itemize}

\medskip

\noindent Finally, let
\begin{equation}
\label{eqn:subalg intro}
\CS^\pm \subset \CV^\pm
\end{equation}
denote the set of elements which pair trivially with $K^\pm$ under the pairings \eqref{eqn:pair intro 1} and \eqref{eqn:pair intro 2}, respectively. The following is our first main result.

\medskip

\begin{theorem}
\label{thm:shuffle}

(Theorem \ref{thm:equal}) We have $\CS^\pm = \oCS^\pm$, thus \eqref{eqn:upsilon intro} yields isomorphisms
\begin{equation}
\label{eqn:iso intro}
\Upsilon^\pm : \UUpm \xrightarrow{\sim} \CS^\pm
\end{equation}
Moreover, the pairings \eqref{eqn:pair intro 1}--\eqref{eqn:pair intro 2} descend to pairings
\begin{align}
&\UUp \otimes \CS^- \xrightarrow{\langle \cdot, \cdot \rangle} \BK \label{eqn:pair 1} \\
&\CS^+ \otimes \UUm \xrightarrow{\langle \cdot, \cdot \rangle} \BK \label{eqn:pair 2}
\end{align} 
which are non-degenerate in both arguments, and coincide under			 \eqref{eqn:iso intro}.

\end{theorem}

\medskip

\noindent If the rational functions \eqref{eqn:general zeta intro} have the property that
\begin{equation}
\label{eqn:symmetric intro}
\frac {\zeta_{ij}(x)}{\zeta_{ji}(x^{-1})}
\end{equation}
is regular and non-zero at $x = \infty$ for all $i,j \in I$, then we recall in Subsection \ref{sub:double} the usual construction of gluing the halves $\UUpm$ into a quantum loop group
\begin{equation}
\label{eqn:double quantum group intro}
\UU = \UUp \otimes \UUo \otimes \UUm
\end{equation}
using the Drinfeld double construction (above, $\UUo$ is a commutative subalgebra generated by symbols $\{h_{i,d}, h'_{i,-d}\}_{i \in I, d\geq 0}$, modulo the relation $h_{i,0}h'_{i,0} = 1, \forall i \in I$). 

\medskip

\subsection{}
 
Before we move on to the other main results of the present paper, let us recall two important particular cases of Theorem \ref{thm:shuffle}, which have already appeared in the literature. Let $I$ be a finite set, assumed to be totally ordered in Definition \ref{def:loop intro}.

\medskip

\begin{definition}
\label{def:loop intro}

Consider $\BK = \BQ(q)$ and 
\begin{equation}
\label{eqn:zeta loop intro}
\zeta_{ij}(x) = \frac {(q^{-d_{ij}} - x) (-x)^{-\delta_{i > j}}}{(1-x)^{\delta_{ij}}}
\end{equation}
where $\{d_{ij}\}$ is the Cartan matrix of a simply-laced \footnote{Recall that this means that $\{d_{ij}\}_{i,j\in I}$ are integers such that $$\begin{cases} d_{ij} = d_{ji} \leq 0 &\text{if }i \neq j \\ d_{ii} = 2 &\text{for all }i \in I \end{cases}$$} Kac-Moody Lie algebra $\fg$. Then
$$
\UU = U_q(L\fg)
$$
is the quantum loop group associated to $\fg$ (see \cite{Loop} for details).

\end{definition} 
 
\medskip

\noindent When $U_q(L\fg)$ is replaced by the usual quantum group $U_q(\fg)$, constructions analogous to those of Subsection \ref{sub:general intro} were performed in \cite{G,R,S}, although they are technically substantially different from those in the present paper.

\medskip

\begin{definition}
\label{def:preprojective intro}

Let $Q$ be a quiver with vertex set $I$ and edge set $E$. Consider $\BK = \BQ(q,t_e)_{e \in E}$ and 
\begin{equation}
\label{eqn:zeta preprojective intro}
\zeta_{ij}(x) = \left(\frac {1-xq^{-1}}{1-x} \right)^{\delta_{ij}} \prod_{e = \oij} \left(\frac 1{t_e}-x \right) \prod_{e = \oji} \left(1 - \frac {t_e}{xq} \right)
\end{equation}
Then it was shown in \cite{Quiver 1} that $\UUp$ is isomorphic to the localized preprojective $K$-theoretic Hall algebra of $Q$ (defined as in \cite{SV Hilb}).

\end{definition}
 
\medskip 

\noindent When $Q$ is the quiver with one vertex and $g$ loops, the algebra $\UUp$ from Definition \ref{def:preprojective intro} also matches the Hall algebra of a genus $g$ curve over the finite field $\BF_{q^{-1}}$ (\cite{Quiver 3}).

\medskip

\subsection{}
\label{sub:duality}

To complete the picture of Subsection \ref{sub:general intro}, one would like an explicit presentation of the defining relations of $\UUpm$, i.e. a collection of generators for the two-sided ideals $K^\pm$. In the setting of Definition \ref{def:loop intro} for $\fg$ of finite type, such a system of generators is given by the Drinfeld-Serre relations for all $i \neq j$ in $I$. However, for arbitrary symmetric Cartan matrices, one needs the more general expressions constructed in \cite[formula (1.15)]{Loop} for any $i \neq j$ in $I$ and integers $k,l \geq 0$ such that
\begin{equation}
\label{eqn:indexing loop}
k+l = -d_{ij}
\end{equation}
Similarly, in the setting of Definition \ref{def:preprojective intro}, the ideals $K^\pm$ are generated by the cubic expressions defined in \cite[formula (1.6)]{Quiver 3} for every arrow of the doubled quiver. Rather than present the aforementioned expressions explicitly, we will invoke the following principle that follows from the non-degeneracy of the pairings \eqref{eqn:pair intro 1}--\eqref{eqn:pair intro 2}
\begin{multline}
\label{eqn:principle}
\textbf{elements of the ideals }K^\pm \subset \tUUpm \textbf{ are dual to} \\ \textbf{linear conditions cutting out the inclusions } \CS^\mp \subset \CV^\mp
\end{multline}
With this in mind, we will find it more informative to present the latter linear conditions. In the setting of Definition \ref{def:loop intro}, for any $i \neq j$ in $I$ and $k,l$ as in \eqref{eqn:indexing loop}, Laurent polynomials $R \in \CS^\mp$ must satisfy the condition
\begin{equation}
\label{eqn:vanish loop}
R \Big|^{z_{i1} = xq^k, z_{i2} = xq^{k-2}, \dots, z_{i,k+1} = xq^{-k}}_{z_{j1} = xq^l, z_{j2} = xq^{l-2}, \dots, z_{j,l+1} = xq^{-l}} = 0
\end{equation}
In the setting of Definition \ref{def:preprojective intro}, Laurent polynomials $R \in \CS^\mp$ satisfy the condition
\begin{equation}
\label{eqn:vanish preprojective}
R\Big|_{z_{ia} = \frac {qz_{jb}}{t_e} = q z_{ic}} = R\Big|_{z_{ja} = t_e z_{ib} = q z_{jc}} = 0
\end{equation}
for any arrow $e$ from $i$ to $j$, and for any $a \neq c$ (and moreover $a\neq b \neq c$ if $i=j$). Comparing \eqref{eqn:vanish loop} to \eqref{eqn:vanish preprojective} makes it apparent that the linear conditions that cut out $\CS^\mp \subset \CV^\mp$ (and dually, the generators of the ideals $K^\pm$) strongly depend on the particular $\zeta$ functions that define our quantum loop groups and shuffle algebras.
 
\medskip

\subsection{} 
\label{sub:main theorem shuffle}

The common feature of the rational functions \eqref{eqn:zeta loop intro} and \eqref{eqn:zeta preprojective intro} is that their numerators are completely split. Therefore, we will henceforth assume that
\begin{equation}
\label{eqn:factor zeta intro}
\zeta_{ij}(x) = \frac {\alpha_{ij} x^{s_{ij}}}{(1-x)^{\delta_{ij}}}\prod_{e=1}^{\#_{ij}} (1-x q_e^{ij})
\end{equation}
for certain $\alpha_{ij}, q_e^{ij} \in \BK^\times$, integers $s_{ij}$ and non-negative integers $\#_{ij}$ \footnote{This is not a major restriction, as completely general rational functions $\zeta_{ij}(x)$ can be written as in \eqref{eqn:factor zeta intro}, at the cost of replacing the ground field $\BK$ by an algebraic field extension.}. Recall from Subsection \ref{sub:duality} that in the special cases studied in Definitions \ref{def:loop intro} and \ref{def:preprojective intro}, the inclusion $\CS^\pm \subset \CV^\pm$ can be presented as the intersection of ideals of Laurent polynomials satisfying conditions  such as \eqref{eqn:vanish loop} and \eqref{eqn:vanish preprojective}, respectively. Our next main Theorem shows that this kind of behavior holds in general. For any $\bn \in \nn$, we will write $\CV_{\pm \bn}$ for the $\bn$-th direct summand of \eqref{eqn:big shuf intro}, and set $\CS_{\pm \bn} = \CS^\pm \cap \CV_{\pm \bn}$. 

\medskip

\begin{theorem}
\label{thm:roots 1}

For any $\bn = (n_i)_{i \in I} \in \nn$ with $n = \sum_{i \in I} n_i$, and any
\begin{equation}
\label{eqn:point intro}
p = (p_{ia})_{i \in I, a \in \{1,\dots, n_i\}} \in (\BK^{\times})^n \Big / \BK^\times
\end{equation}
there exists a homogeneous ideal
\begin{equation}
\label{eqn:ideal intro}
\CI_{p, \bn} \subseteq \BK[z_{i1}, z_{i1}^{- 1}, \dots, z_{in_i}, z_{in_i}^{- 1}]_{i \in I}
\end{equation}
supported on the one-dimensional locus $(z_{ia})_{i \in I, a \in \{1,\dots, n_i\}} \in p\BK^\times$, such that
\begin{equation}
\label{eqn:shuffle intersect intro}
\CS_{\pm \bn} = \CV_{\pm \bn} \mathop{\bigcap_{\b0 \leq \bk \leq \bn}}_{p \in (\BK^{\times})^k / \BK^\times} \CI_{p, \bk}^{(\bn)}
\end{equation}
(see \eqref{eqn:big int} and the discussion preceding it for the notation in the right-hand side). 

\end{theorem}

\medskip

\subsection{} 
\label{sub:main theorem quantum}

According to the principle \eqref{eqn:principle}, dualizing Theorem \ref{thm:roots 1} will yield a collection of generators for the ideals $K^\pm$ (and thus complete the generators-and-relations presentation of the quantum loop groups $\UU$). Let us  write
$$
\tUUpm = \bigoplus_{\bn \in \nn} \tUU_{\pm \bn} = \bigoplus_{(\bn,d) \in \nn \times \BZ} \tUU_{\pm \bn, \pm d}
$$
for the components with respect to the grading defined by $\deg e_{i,d} = (\bs^i,d)$, $\deg f_{i,d} = (-\bs^i,d)$ for all $i \in I$ and $d \in \BZ$, where $\bs^i = \underbrace{(0,\dots,0,1,0,\dots,0)}_{1\text{ on }i\text{-th position}} \in \nn$.

\medskip

\begin{theorem}
\label{thm:roots 2}

For any $(\bn,d) \in \nn \times \BZ$ and $p \in (\BK^{\times})^n / \BK^\times$, there are finite sets
\begin{equation}
\label{eqn:elements intro}
W_{p,\pm \bn, \pm d} \subset \tUU_{\pm \bn,\pm d}
\end{equation}
which generate the two-sided ideal $K^\pm$ as $p$ and $(\bn,d)$ vary. Thus, we have
\begin{equation}
\label{eqn:quantum intro define}
\UUpm = \tUUpm \Big/ \left(W_{p,\pm \bn, \pm d}\right)_{p \in (\BK^{\times})^n / \BK^\times, \bn \in \nn, d\in \BZ}
\end{equation}


\end{theorem}

\medskip

\noindent As we will show in the proof of Proposition \ref{prop:primary}, one can compute the ideals $\CI_{p, \bn}$ recursively in $\bn$. In more detail, we emphasize the fact that each ideal \eqref{eqn:ideal intro} is homogeneous and supported on the one-dimensional locus $p \BK^\times \cong \text{Spec } \BK[x^{\pm 1}]$. Therefore, in every homogeneous degree $d \in \BZ$, the inclusion
$$
\CI_{p, \bn, d} \subseteq \BK[z_{i1},z_{i1}^{-1} \dots, z_{in_i}, z_{in_i}^{- 1}]_{i \in I}^{\text{homogeneous degree d}}
$$
is cut out as a $\BK$-vector space by finitely many linear conditions. On $\CV_{\pm \bn}$, these linear conditions are realized by pairing with finitely many elements of $\tUU_{\mp \bn,\mp d}$ under the pairings \eqref{eqn:pair intro 1}--\eqref{eqn:pair intro 2}, as shown in the proof of Proposition \ref{prop:dual}. The aforementioned elements are, by definition, the elements of the finite sets \eqref{eqn:elements intro}.

\medskip

\subsection{} 
\label{sub:wheel}

Because of the arbitrariness of the point \eqref{eqn:point intro}, it might seem like one must consider a great multitude of ideals \eqref{eqn:ideal intro} and elements \eqref{eqn:elements intro}. However, this is not the case: as we will see in Subsection \ref{sub:point}, for every $\bn \in \nn$ and every $p = (p_{ia})_{i \in I, a \in \{1,\dots,n_i\}}$ the corresponding ideals/elements only depend on the collection 
$$
\left\{q_e^{ij} \Big| q_e^{ij} = \frac {p_{jb}}{p_{ia}} \text{ for some }a,b \right\}_{i,j \in I, e \in \{1,\dots,\#_{ij}\}} \subset \Big\{ \text{roots of the }\zeta_{ij}\text{'s} \Big\}
$$
As there are only finitely many choices of such collections, one only needs to consider finitely many ideals/elements in \eqref{eqn:ideal intro}/\eqref{eqn:elements intro} for every $\bn \in \nn$. Moreover, the ideal \eqref{eqn:ideal intro} is non-trivial (respectively the set \eqref{eqn:elements intro} is non-empty) only if the point $p$ is a wheel, i.e. there exists a cycle $(i_1,a_1), \dots, (i_k,a_k), (i_{k+1},a_{k+1}) = (i_1,a_1)$ s.t.
\begin{equation}
\label{eqn:wheel}
\frac {p_{i_{\bullet + 1} a_{\bullet + 1}}}{p_{i_\bullet a_{\bullet}}} \in \Big\{q_e^{i_\bullet i_{\bullet+1}} \Big\}_{e \in \{1,\dots,\#_{i_\bullet i_{\bullet+1}}\}}, \qquad \forall \ \bullet \in \{1,\dots,k\}
\end{equation}
For such a wheel, the requirement that a Laurent polynomial $R(z_{ia})_{i \in I, a \in \{1,\dots,n_i\}}$ lies in $\CI_{p,\bn}$ is a linear condition on the derivatives of $R$ at the points of the form
$$
(z_{ia})_{i \in I, a \in \{1,\dots,n_i\}} \in p\BK^\times
$$
thus generalizing the \textbf{wheel conditions} discovered in \cite{FO} in the setting of quantum loop groups of finite and affine type (which involved certain special points $p$). In special cases (such as \eqref{eqn:vanish loop}, \eqref{eqn:vanish preprojective} or \eqref{eqn:vanish qqta} below), we see that wheel conditions are rather simple vanishing conditions, but in general they can be quite complicated.

\medskip

\subsection{}
\label{sub:k-ha intro}

Although we do not explicitly describe the wheel conditions for arbitrary $\zeta$ functions, in Subsection \ref{sub:point} we will explain that they are built out of the special cases when
\begin{equation}
\label{eqn:k-ha intro}
\zeta_{ij}(x) = (1-x)^{\#_{ij}-\delta_{ij}}
\end{equation}
for various non-negative integers $\#_{ij}$. In this setting, the shuffle algebra $\CV^+$ is none other than the $K$-theoretic Hall algebra with 0 potential associated to the quiver $Q$ with vertex set $I$ and $\#_{ij}$ arrows from the vertex $i$ to the vertex $j$, for all $i,j \in I$ (see \cite{KS, Pad 1}). These algebras are very actively studied in connection with Donaldson-Thomas invariants for Calabi-Yau categories, and the upshot of our theorems is that the ``spherical" subalgebra $\oCS^+ = \CS^+ \subset \CV^+$ admits a presentation as (half of) a quantum loop group. 

\medskip

\noindent As for $K$-theoretic Hall algebras equivariant with respect the torus 
$$
T = \prod_{i,j \in I} (\BC^*)^{\#_{ij}}
$$
we let $q^{ij}_{e}$ denote the standard character corresponding to the $e$-th factor of $\BC^*$ in the $i,j$ term of the product above. Then let us work over the ring
$$
\BL = \BQ \left[ \left(q_e^{ij} \right)^{\pm 1} \right]_{i,j \in I, e \in \{1,\dots,\#_{ij}\}}
$$
and define the $\BL$-algebra
$$
\CV_{\text{int}}^+ = \bigoplus_{\bn = (n_i)_{i \in I} \in \nn} \BL [z_{i1}, z_{i1}^{-1}, \dots, z_{in_i}, z_{in_i}^{-1}]_{i \in I}^{\sym}
$$
endowed with the shuffle multiplication \eqref{eqn:shuf prod} associated to the rational functions \eqref{eqn:factor zeta intro} \footnote{One can define $\alpha_{ij}$ and $s_{ij}$ arbitrarily in \eqref{eqn:factor zeta intro}, and various choices correspond to various twists of the $K$-theoretic Hall product; the customary choices are $\alpha_{ij} = 1$ and $s_{ij} = 0$ for all $i,j$.}. Similarly, one lets 
\begin{equation}
\label{eqn:integral shuffle}
\CS_{\text{int}}^+ \subseteq \CV_{\text{int}}^+
\end{equation}
denote the $\BL$-subalgebra generated by $\{z_{i1}^d\}_{i \in I, d \in \BZ}$. The methods in the present paper allow to describe the fibers of the subalgebra \eqref{eqn:integral shuffle} above various $\overline{\BQ}$-points of $\text{Spec }\BL$, which corresponds to various specializations of the arrow parameters $q_e^{ij}$. If the specialization is generic, there are no wheels \eqref{eqn:wheel} and thus
$$
\CS^+_{\text{generic}} = \CV^+_{\text{generic}}
$$
(dually, the quantum loop group associated to generic parameters is $\UUm = \tUUm$). At the opposite extreme, if the arrow parameters are all specialized to $1$, we recover the situation of \eqref{eqn:k-ha intro}, which as we will see is the most interesting and important one.

\medskip

\subsection{}
\label{sub:applications intro 2}

\noindent Our construction also applies to quiver quantum toroidal algebras. These are trigonometric versions (introduced in \cite{GLY,NW}) of the quiver Yangians (introduced in \cite{LY}, see also \cite{RSYZ} for a related mathematical construction) that act on the vector spaces of BPS states for non-compact toric Calabi-Yau threefolds $X$. More specifically, one associates to such an $X$ a quiver $Q$ endowed with arrow parameters
$$
\{ q_e^{ij} \}_{i,j \in I, e \in \{1,\dots,\#_{ij}\}} \in \text{Rep}_{\BC^* \times \BC^*}
$$
See \cite[Section B]{NW} for the precise construction; the two-dimensional torus $\BC^* \times \BC^*$ should be interpreted as the kernel of the Calabi-Yau form. When the quiver is symmetric (i.e. $\#_{ij} = \#_{ji}$, which is known as ``non-chiral" in the physics literature), the corresponding \textbf{quiver quantum toroidal algebra} is none other than our
$$
\tUU = \tUUp \otimes \tUUo \otimes \tUUm
$$
defined with respect to the $\zeta$ functions \eqref{eqn:factor zeta intro}, for certain $\alpha_{ij}$ and $s_{ij}$ \footnote{For non-symmetric (i.e. chiral) quivers, our construction describes the positive/negative parts of quiver quantum toroidal algebras, but the Hopf algebra structure is not well-defined.}. One of the main features of these quiver quantum toroidal algebras is that they act on the vector spaces of BPS crystal configurations 
$$
M = \bigoplus_{\Lambda \text{ 3d crystal configuration}} \BK \cdot |\Lambda\rangle  
$$
where $\BK = \text{Frac}(\text{Rep}_{\BC^* \times \BC^*})$ is the ground field in the present setup (see \cite[Section 5]{NW} for a review of 3d crystal configurations, which are generalizations of plane partitions). As noted in \cite[Section 5]{GLY}, the action
$$
\tUUpm \curvearrowright M
$$
factors through the homomorphism $\tUUpm \xrightarrow{\tUpsilon^\pm} \CS^\pm \stackrel{\eqref{eqn:iso intro}}\cong \UUpm$. Thus, there is an action
$$
\UU \curvearrowright M
$$
where $\UU$ is the Drinfeld double \eqref{eqn:double quantum group intro}. In other words, all the elements \eqref{eqn:elements intro} act by 0 in the representation $M$, and so the quotient $\UU$ plays the role of a ``reduced" version of the quiver quantum toroidal algebra $\tUU$. In \cite[formula (1.13)]{QQTA}, we will give a generators-and-relations presentation of the quantum loop group $\UU$ in the setting at hand, by giving an explicit set of generators for the ideals $K^\pm$. For now let us mention the dual statement, in accordance with the principle \eqref{eqn:principle}: a Laurent polynomial $R \in \CV^\mp$ lies in $\CS^{\mp}$ if and only if satisfies the conditions
\begin{equation}
\label{eqn:vanish qqta}
R\Big|_{z_a = z_{a-1} q_{e_a}^{i_ai_{a-1}}, \ \forall a \in \{1,\dots,k\}} = 0
\end{equation}
for any face $F = \{i_0,i_1,i_2,\dots,i_{k-1},i_k=i_0\}$ of the quiver $Q$ (in the context of quiver quantum toroidal algebras, the quiver $Q$ is naturally drawn on the torus with polygonal faces), whose boundary edges are denoted by $e_1,\dots,e_k$.

\medskip

\subsection{}
\label{sub:plan}

The outline of the present paper is the following.

\medskip

\begin{itemize}[leftmargin=*] 

\item In Section \ref{sec:shuffle}, we prove Theorem \ref{thm:shuffle} for arbitrary $\zeta$ functions \eqref{eqn:general zeta intro} 

\medskip

\item In Section \ref{sec:roots}, we prove Theorems \ref{thm:roots 1} and \ref{thm:roots 2} for factored $\zeta$ functions, as in \eqref{eqn:factor zeta intro}

\end{itemize}

\medskip

\noindent I would like to express my gratitude to all the mathematicians and physicists that have worked on the wonderful constructions and ideas referenced in the present paper, with special thanks to Igor Frenkel for his many beautiful and fundamental contributions. I gratefully acknowledge NSF grant DMS-$1845034$, as well as support from the Alfred P.\ Sloan Foundation and the MIT Research Support Committee.

\bigskip

\section{Shuffle algebras} 
\label{sec:shuffle}

\medskip

\noindent We will develop the basic theory of shuffle algebras, in the trigonometric setting studied by \cite{E, FO}, which depend on a choice of rational functions $\zeta_{ij}(x)$. In the present Section we will deal mostly with generalities that apply to arbitrary $\zeta_{ij}$'s.

\medskip 

\subsection{} 
\label{sub:def quad}

We will work over a base field $\BK$ of characteristic 0. Consider a finite set $I$, and let us fix a non-zero Laurent polynomial for any $i \neq j$ in $I$
\begin{equation}
\label{eqn:def zeta diff}
\zeta_{ij}(x) \in \BK[x^{\pm 1}]
\end{equation}
and a non-zero rational function with at most a simple pole at $x = 1$ for any $i \in I$
\begin{equation}
\label{eqn:def zeta same}
\zeta_{ii}(x) \in \frac {\BK[x^{\pm 1}]}{1-x}
\end{equation}

\medskip

\begin{definition}
\label{def:quad}

The (positive part of the) \textbf{quadratic quantum loop group} $\tUUp$ associated to the datum $\{\zeta_{ij}\}_{i,j \in I}$ is the $\BK$-algebra generated by symbols
$$
\{e_{i,d}\}_{i \in I, d \in \BZ}
$$
modulo the following relations for all $i,j \in I$
\begin{equation}
\label{eqn:rel quad}
e_i(z) e_j(w) \zeta_{ji} \left(\frac wz\right) = e_j(w) e_i(z) \zeta_{ij} \left( \frac zw \right)
\end{equation}
Above and henceforth, we use the notation:
$$
e_i(z) = \sum_{d \in \BZ} \frac {e_{i,d}}{z^d}
$$
for all $i \in I$, and relation \eqref{eqn:rel quad} is interpreted as an infinite collection of relations obtained by equating the coefficients of all $\{z^aw^b\}_{a,b\in \BZ}$ in the left and right-hand sides (if $i = j$, one clears the denominators $z-w$ from \eqref{eqn:rel quad} before equating coefficients). 

\end{definition}

\medskip

\noindent The algebra $\tUUp$ is graded by $\nn \times \BZ$ (in this paper, $\BN$ is assumed to contain 0), via
$$
\deg e_{i,d} = (\bs^i, d)
$$
for all $i \in I$ and $d \in \BZ$. Above and throughout the present paper, $\bs^i \in \nn$ denotes the $I$-tuple of integers with a 1 on position $i$, and 0 everywhere else. We will write
$$
\tUUp = \bigoplus_{\bn \in \nn} \tUU_{\bn} = \bigoplus_{(\bn,d) \in \nn \times \BZ} \tUU_{\bn,d}
$$
for the graded components. We have the \textbf{shift automorphism}
\begin{equation}
\label{eqn:shift 1}
\tUUp \xrightarrow{\tau_{\bk}} \tUUp, \qquad e_{i,d} \mapsto e_{i,d+k_i}
\end{equation}
for any $\bk = (k_i)_{i \in I} \in \zz$.

\medskip

\subsection{}
\label{sub:def shuf}

Let us consider an infinite collection of variables $z_{i1},z_{i2},\dots$ for all $i \in I$. For any $\bn = (n_i)_{i \in I} \in \nn$, we will write $|\bn| = \sum_{i \in I} n_i$ and $\bn! = \prod_{i \in I} n_i!$. The following construction is a straightforward generalization of that of \cite{E, FO}.

\medskip

\begin{definition}
\label{def:big shuf}

The \textbf{big shuffle algebra} associated to the datum $\{\zeta_{ij}\}_{i,j \in I}$ is
$$
\CV^+ = \bigoplus_{\bn \in \BN^I} \BK[z_{i1}^{\pm 1},\dots,z_{in_i}^{\pm 1}]_{i \in I}^{\esym}
$$
endowed with the multiplication
\begin{equation}
\label{eqn:shuf prod}
R(\dots,z_{i1},\dots,z_{in_i},\dots) * R'(\dots,z_{i1},\dots,z_{in_i'},\dots ) = 
\end{equation}
$$
\eSym \left[ \frac {R(\dots,z_{i1},\dots,z_{in_i},\dots) R'(\dots,z_{i,n_i+1},\dots,z_{i,n_i+n_i'},\dots)}{\bn! \bn'!} \mathop{\prod^{i,j \in I}_{1\leq a\leq n_i}}_{n_j < b \leq n_j+n_j'} \zeta_{ij} \left(\frac {z_{ia}}{z_{jb}} \right) \right]
$$
Above and henceforth, ``\emph{sym}" (resp. ``\emph{Sym}") denotes symmetric functions (resp. symmetrization) with respect to the variables $z_{i1},z_{i2},\dots$ for each $i \in I$ separately \footnote{Although the $\zeta$ functions might seem to contribute simple poles at $z_{ia}-z_{ib}$ for $a \neq b$ to the right-hand side of \eqref{eqn:shuf prod}, these poles disappear when taking the symmetrization (the poles in question can only have even order in any symmetric rational function).}. 

\end{definition} 

\medskip

\noindent Note that the algebra $\CV^+$ is graded by $\nn \times \BZ$, via
$$
\deg R(\dots, z_{i1}, \dots, z_{in_i}, \dots) = (\bn, \text{hom deg }R)
$$
where ``$\text{hom deg }R$" denotes the homogeneous degree of $R$ in all its variables. Let
$$
\CV^+ = \bigoplus_{\bn \in \nn} \CV_{\bn} = \bigoplus_{(\bn,d) \in \nn \times \BZ} \CV_{\bn,d}
$$
denote the graded pieces. The analogue of the shift automorphism \eqref{eqn:shift 1} is
\begin{equation}
\label{eqn:shift 2}
\CV^+ \xrightarrow{\tau_{\bk}} \CV^+, \qquad R(\dots, z_{ia}, \dots) \mapsto R(\dots, z_{ia}, \dots) \prod_{i \in I, a \geq 1} z_{ia}^{k_i}
\end{equation}
for any $\bk = (k_i)_{i \in I} \in \zz$.

\medskip

\subsection{} We will also encounter the opposite algebra
$$
\tUUm = \BK \Big \langle f_{i,d} \Big \rangle_{i \in I, d \in \BZ} \Big / \text{relation \eqref{eqn:opposite}}
$$
with the defining relation
\begin{equation}
\label{eqn:opposite}
f_i(z)f_j(w) \zeta_{ij} \left( \frac zw \right) = f_j(w) f_i(z) \zeta_{ji}\left(\frac wz \right)
\end{equation}
for all $i,j \in I$. Also let $\CV^- = \CV^+$ as a vector space, but in the multiplication \eqref{eqn:shuf prod} we replace
$$
\zeta_{ij} \left(\frac {z_{ia}}{z_{jb}} \right) \quad \text{with} \quad \zeta_{ji} \left(\frac {z_{jb}}{z_{ia}} \right)
$$
The algebra $\tUUm$ is $(-\nn) \times \BZ$ graded, via
$$
\deg f_{i,d} = (-\bs^i, d)
$$
for all $i \in I$, $d \in \BZ$, while $\CV^-$ is $(-\nn) \times \BZ$ graded, via
$$
\deg R(\dots, z_{i1}, \dots, z_{in_i}, \dots) = (-\bn, \text{hom deg }R)
$$
for all $R \in \CV^-$. We will denote the graded pieces of $\tUUm$ and $\CV^-$ by
\begin{align*}
&\tUUm = \bigoplus_{\bn \in \nn} \tUU_{-\bn} = \bigoplus_{(\bn, d) \in \nn \times \BZ} \tUU_{-\bn,d} \\
&\CV^- = \bigoplus_{\bn \in \nn} \CV_{-\bn} = \bigoplus_{(\bn, d) \in \nn \times \BZ} \CV_{-\bn,d}
\end{align*}
The shift automorphisms on the negative algebras are defined as
\begin{align}
&\tUUm \xrightarrow{\tau_{\bk}} \tUUm, \qquad e_{i,d} \mapsto e_{i,d-k_i} \label{eqn:shift 3} \\
&\CV^- \xrightarrow{\tau_{\bk}} \CV^-, \qquad R(\dots, z_{ia}, \dots) \mapsto R(\dots, z_{ia}, \dots) \prod_{i \in I, a \geq 1} z_{ia}^{-k_i} \label{eqn:shift 4}
\end{align}
for any $\bk = (k_i)_{i \in I} \in \zz$.

\medskip

\subsection{} 
\label{sub:upsilon}

There are homomorphisms of $(\pm \nn) \times \BZ$ graded $\BK$-algebras
\begin{equation}
\label{eqn:tilde upsilon}
\tUUpm \xrightarrow{\tUpsilon^\pm} \CV^\pm, \qquad e_{i,d}, f_{i,d} \mapsto z_{i1}^d
\end{equation}
Indeed, this claim only entails checking the fact that relations \eqref{eqn:rel quad} are respected by the shuffle product \eqref{eqn:shuf prod}, which is straightforward. The map \eqref{eqn:tilde upsilon} is neither injective nor surjective, and one of the main goals of the present paper is to describe
$$
\oCS^\pm = \text{Im }\widetilde{\Upsilon}^\pm
$$
(i.e. $\oCS^\pm$ is the $\BK$-subalgebra of $\CV^\pm$ generated by $\{z_{i1}^d\}_{i \in I, d\in \BZ}$) and
$$
K^\pm = \text{Ker }\widetilde{\Upsilon}^\pm
$$
as a two-sided ideal of $\tUUpm$. 

\medskip

\subsection{} A key role in our study of the algebras $\tUUpm$ and $\CV^\pm$ is played by certain bilinear pairings that we will introduce shortly. Let us consider the following notation for all rational functions $F(z_1,\dots,z_n)$. If $Dz_a = \frac {dz_a}{2\pi i z_a}$, then we will write
\begin{equation}
\label{eqn:contour integral}
\int_{|z_1| \gg \dots \gg |z_n|} F(z_1,\dots,z_n) \prod_{a=1}^n Dz_a
\end{equation}
for the constant term in the expansion of $F$ as a power series in 
$$
\frac {z_2}{z_1}, \dots, \frac {z_n}{z_{n-1}}
$$
The notation in \eqref{eqn:contour integral} is motivated by the fact that if $\BK = \BC$, one could compute this constant term as a contour integral (with the contours being concentric circles, situated very far from each other compared to the absolute values of the coefficients of $F$). We define $\int_{|z_1| \ll \dots \ll |z_n|} F(z_1,\dots,z_n) \prod_{a=1}^n Dz_a$ by analogy with \eqref{eqn:contour integral}.

\medskip

\begin{definition}
\label{def:pair}

There exist bilinear pairings
\begin{align}
&\tUUp \otimes \CV^- \xrightarrow{\langle \cdot, \cdot \rangle} \BK \label{eqn:pair} \\
&\CV^+ \otimes \tUUm \xrightarrow{\langle \cdot, \cdot \rangle} \BK \label{eqn:pair opposite}
\end{align}
given for all $R^\pm \in \CV_{\pm \bn}$ and all $i_1,\dots,i_n \in I$, $d_1,\dots,d_n \in \BZ$ by
\begin{align}
&\Big \langle e_{i_1,d_1} \cdots e_{i_n,d_n}, R^- \Big \rangle = \int_{|z_1| \gg \dots \gg |z_n|} \frac {z_1^{d_1}\dots z_n^{d_n} R^-(z_1,\dots,z_n)}{\prod_{1\leq a < b \leq n} \zeta_{i_bi_a} \left(\frac {z_b}{z_a} \right)} \prod_{a=1}^n Dz_a \label{eqn:pair formula} \\
&\Big \langle R^+, f_{i_1,-d_1} \cdots  f_{i_n,-d_n} \Big \rangle = \int_{|z_1| \ll \dots \ll |z_n|} \frac {z_1^{-d_1}\dots z_n^{-d_n} R^+(z_1,\dots,z_n)}{\prod_{1\leq a < b \leq n} \zeta_{i_ai_b} \left(\frac {z_a}{z_b} \right)} \prod_{a=1}^n Dz_a 
\label{eqn:pair formula opposite}
\end{align}
if $\bs^{i_1}+\dots +\bs^{i_n} = \bn$, and 0 otherwise. 

\end{definition}

\medskip

\noindent In the right-hand sides of \eqref{eqn:pair formula} and \eqref{eqn:pair formula opposite}, we implicitly identify
\begin{equation}
\label{eqn:relabeling}
z_a \quad \text{with} \quad z_{i_a\bullet_a}, \quad \forall a \in \{1,\dots, n\}
\end{equation}
where $\bullet_a \in \{1,2,\dots,n_{i_a}\}$ may be chosen arbitrarily (however, we require $\bullet_a \neq \bullet_b$ if $a\neq b$ and $i_a = i_b$) due to the symmetry of $R^\pm$. We call \eqref{eqn:relabeling} a \textbf{relabeling} of the variables of $R^\pm(z_{i1},\dots,z_{in_i})_{i \in I}$ in accordance with $i_1,\dots,i_n$.

\medskip

\subsection{}

The pairings \eqref{eqn:pair}--\eqref{eqn:pair opposite} are non-degenerate in the $\CV^\pm$ argument, i.e.
\begin{align}
&\Big \langle \tUUp, R^- \Big \rangle = 0 \quad \Rightarrow \quad R^-  = 0 \label{eqn:non-degenerate 1} \\
&\Big \langle R^+,\tUUm \Big \rangle = 0 \quad \Rightarrow \quad R^+  = 0 \label{eqn:non-degenerate 2}
\end{align}
for any $R^\pm \in \CV^\pm$. This is simply because any rational function whose power series expansion (with respect to any order of variables) vanishes must be identically 0. 

\medskip

\begin{definition}
\label{def:shuf}

Let $\CS^\mp \subset \CV^\mp$ denote the dual of $K^\pm = \emph{Ker }\tUpsilon^\pm$ under the pairings \eqref{eqn:pair}--\eqref{eqn:pair opposite}, respectively, i.e.
\begin{align}
&R^- \in \CS^- \quad \Leftrightarrow \quad \Big \langle K^+, R^- \Big \rangle = 0 \label{eqn:shuf 2} \\
&R^+ \in \CS^+ \quad \Leftrightarrow \quad \Big \langle R^+, K^- \Big \rangle = 0 \label{eqn:shuf 1} 
\end{align}
We will call $\CS^\pm$ the positive/negative \textbf{shuffle algebra} \footnote{The fact that $\CS^\pm$ is closed under shuffle product will be established as part of Theorem \ref{thm:equal}.}.

\end{definition}

\medskip

\noindent As a consequence of Definition \ref{def:shuf}, the pairings \eqref{eqn:pair}--\eqref{eqn:pair opposite} descend to pairings
\begin{align}
&\oCS^+ \otimes \CS^- \xrightarrow{\langle \cdot, \cdot \rangle} \BK \label{eqn:descended pairing} \\
&\CS^+ \otimes \oCS^- \xrightarrow{\langle \cdot, \cdot \rangle} \BK \label{eqn:descended pairing opposite}
\end{align}
which are non-degenerate in the $\CS^\pm$ argument. The main result of the present Section is the following.

\medskip

\begin{theorem} 
\label{thm:equal}

We have $\CS^\pm = \oCS^\pm$, and the pairings \eqref{eqn:descended pairing} and \eqref{eqn:descended pairing opposite} coincide. 

\end{theorem}

\medskip

\subsection{} Before we delve into proving Theorem \ref{thm:equal}, let us note that it states that
$$
R \in \oCS^- \qquad \Leftrightarrow \qquad R \in \CS^-
$$
for an arbitrary $R \in \CV^-$. The first $\in$ above is an issue of the Laurent polynomial $R$ lying in a certain ideal, while the second $\in$ above involves checking that $R$ is annihilated by a collection of linear maps. Explicitly, if $\bn \in \nn$ denotes the (negative of the) degree of $R$, let us choose a linear spanning set
\begin{equation}
\label{eqn:generator explicit}
\mathop{\sum_{i_1,\dots,i_n \in I}}_{\bs^{i_1}+\dots+\bs^{i_n} = \bn} \Big[ p_{i_1,\dots,i_n}(z_1,\dots,z_n) e_{i_1}(z_1) \dots e_{i_n}(z_n) \Big]_{\text{ct}} \in K^+ \cap \tUU_{\bn}
\end{equation}
(where $p_{i_1,\dots,i_n}$ are Laurent polynomials, and $[\dots]_{\text{ct}}$ denotes the constant term in the variables $z_1,\dots,z_n$). If we relabel the variables according to \eqref{eqn:relabeling}, then \eqref{eqn:generator explicit} is equivalent to the symmetric Laurent polynomial identity
\begin{equation}
\label{eqn:generator explicit equivalent}
\mathop{\sum_{i_1,\dots,i_n \in I}}_{\bs^{i_1}+\dots+\bs^{i_n} = \bn} \Sym  \left[ p_{i_1,\dots,i_n}(z_{i_1\bullet_1},\dots,z_{i_n \bullet_n}) \prod_{1\leq a < b \leq n} \zeta_{i_ai_b} \left( \frac {z_{i_a \bullet_a}}{z_{i_b \bullet_b}} \right) \right] = 0
\end{equation}
Then $R \in \CV_{-\bn}$ lies in $\CS_{-\bn} = \CV_{-\bn} \cap \CS^{-}$ if and only if for any $\{p_{i_1,\dots,i_n}\}_{i_1,\dots,i_n \in I}$ as in \eqref{eqn:generator explicit}, we have
$$
0 = \left \langle \mathop{\sum_{i_1,\dots,i_n \in I}}_{\bs^{i_1}+\dots+\bs^{i_n} = \bn} \Big[ p_{i_1,\dots,i_n}(z_1,\dots,z_n) e_{i_1}(z_1) \dots e_{i_n}(z_n) \Big]_{\text{ct}} , R \right \rangle =
$$
\begin{equation}
\label{eqn:to check vector space}
= \mathop{\sum_{i_1,\dots,i_n \in I}}_{\bs^{i_1}+\dots+\bs^{i_n} = \bn} \int_{|z_1| \gg \dots \gg |z_n|} \frac {p_{i_1,\dots,i_n}(z_1,\dots,z_n) R(z_1,\dots,z_n)}{\prod_{1\leq a < b \leq n} \zeta_{i_bi_a} \left(\frac {z_b}{z_a} \right)} \prod_{a=1}^n Dz_a
\end{equation}
Formula \eqref{eqn:to check vector space} yields a collection of equations, linear in $R$, which determine the $\BK$-vector subspace $\CS_{-\bn} \subset \CV_{-\bn}$ for any given $\bn \in \nn$. The situation of the $\BK$-vector subspace $\CS_{\bn} \subset \CV_{\bn}$ is treated analogously, with $f$'s instead of $e$'s.

\medskip

\begin{example} 
\label{ex:basic}

Assume $n=3$, $i_1,i_2,i_3 \in I$ are distinct and
\begin{align*}
&\zeta_{i_1i_2}(x) = \zeta_{i_2i_3}(x) = \zeta_{i_3i_1}(x) = 1 \\
&\zeta_{i_2i_1}(x) = \zeta_{i_3i_2}(x) = \zeta_{i_1i_3}(x) = 1-x
\end{align*}
Then every coefficient of the formal series
\begin{equation}
\label{eqn:formal series (1,1,1)}
e_{i_1}(z_1)e_{i_2}(z_2)e_{i_3}(z_3) + \frac {z_1}{z_3} \cdot e_{i_2}(z_2)e_{i_3}(z_3) e_{i_1}(z_1) + \frac {z_2}{z_3} \cdot e_{i_3}(z_3) e_{i_1}(z_1)e_{i_2}(z_2)
\end{equation}
lies in $K^+$, and property \eqref{eqn:to check vector space} is equivalent to
\begin{equation}
\label{eqn:pair formal series (1,1,1)}
\delta \left(\frac {z_1}{z_2} \right)\delta \left(\frac {z_1}{z_3} \right)R(z_1,z_2,z_3) = 0 \quad \Leftrightarrow \quad R(x,x,x) = 0
\end{equation}
where $\delta(x) = \sum_{d \in \BZ} x^d$ is the formal delta series.

\end{example}

\medskip

\subsection{} Given $\bn \in \nn$, let us now fix $i_1,\dots,i_n \in I$ such that $\bn = \bs^{i_1}+\dots+\bs^{i_n}$. Let
$$
S(\bn) \subset S(n)
$$
be the subset of permutations $\sigma$ such that $\sigma(a) < \sigma(b)$ if $a < b$ and $i_a = i_b$. In the present Subsection, we will relabel our variables according to \eqref{eqn:relabeling}. It is easy to see that \eqref{eqn:to check vector space} is equivalent to
\begin{equation}
\label{eqn:linear relations}
0 = \sum_{\sigma \in S(\bn)} \int_{|z_{\sigma(1)}| \gg \dots \gg |z_{\sigma(n)}|} \frac {p_{\sigma}(z_{\sigma(1)}, \dots, z_{\sigma(n)}) R (z_1, \dots, z_n)}{\prod_{1\leq a < b \leq n} \zeta_{i_{\sigma(b)}i_{\sigma(a)}} \left(\frac {z_{\sigma(b)}}{z_{\sigma(a)}} \right)} \prod_{a=1}^n Dz_a
\end{equation}
where we choose to index the Laurent polynomial $p$ by an arbitrary permutation $\sigma \in S(\bn)$ which permutes the fixed indices $i_1,\dots,i_n$. In the previous Subsection, we showed that $R \in \CV_{-\bn}$ lies in $\CS_{-\bn}$ if and only if \eqref{eqn:linear relations} holds for any collection of Laurent polynomials $\{p_{\sigma}\}_{\sigma \in S(\bn)}$ that belongs to the kernel of the following map 
\begin{equation}
\label{eqn:phi}
\bigoplus_{\sigma \in S(\bn)} \BK[z_1^{\pm 1}, \dots, z_n^{\pm 1}] \xrightarrow{\Phi}  \BK[z_1^{\pm 1}, \dots, z_n^{\pm 1}]^{\sym}
\end{equation}
$$
\Phi \Big( \{ p_{\sigma} \}_{\sigma \in S(\bn)}\Big) =  \sum_{\sigma \in S(\bn)} \Sym \left[ p_{\sigma}(z_{\sigma(1)},\dots,z_{\sigma(n)}) \prod_{1\leq a < b \leq n} \zeta_{i_{\sigma(a)}i_{\sigma(b)}} \left( \frac {z_{\sigma(a)}}{z_{\sigma(b)}} \right) \right]
$$
of $\BK[z_1^{\pm 1}, \dots, z_n^{\pm 1}]^{\sym}$-modules (in the present Subsection, ``sym" and ``Sym" refer to symmetric polynomials with respect to any variables $z_a$ and $z_b$ for which $i_a = i_b$).

\medskip

\begin{proposition}
\label{prop:trivial}

If $\{p_{\sigma}\}_{\sigma \in S(\bn)}$ lies in the image of the map 
\begin{equation}
\label{eqn:psi}
\bigoplus_{\sigma \neq \sigma' \in S(\bn)} \BK[z_1^{\pm 1}, \dots, z_n^{\pm 1}] \xrightarrow{\Psi} \bigoplus_{\sigma \in S(\bn)} \BK[z_1^{\pm 1}, \dots, z_n^{\pm 1}]
\end{equation}
$$
\Psi \Big ( \{f_{\sigma\sigma'}(z_1, \dots, z_n) \}_{\sigma \neq \sigma' \in S(\bn)} \Big) =  \{p_{\sigma}(z_1,\dots,z_n)\}_{\sigma \in S(\bn)}
$$
where each $p_{\sigma}(z_1,\dots,z_n)$ is given by
\begin{multline*}
\sum_{\sigma' \neq \sigma} f_{\sigma\sigma'}(z_{\sigma^{-1}(1)}, \dots, z_{\sigma^{-1}(n)}) \prod^{a<b}_{{\sigma'}^{-1} ( \sigma(a) ) > {\sigma'}^{-1} ( \sigma(b) )} \tzeta_{i_bi_a} \left(\frac {z_b}{z_a} \right)  - \\
- \sum_{\sigma' \neq \sigma} f_{\sigma'\sigma}(z_{{\sigma}^{-1}(1)}, \dots, z_{{\sigma}^{-1}(n)}) \prod^{a<b}_{{\sigma'}^{-1} ( \sigma(a) ) > {\sigma'}^{-1} ( \sigma(b) )} \tzeta_{i_bi_a} \left(\frac {z_b}{z_a} \right) \left(- \frac {z_a}{z_b} \right)^{\delta_{i_ai_b}}
\end{multline*}
then the right-hand side of \eqref{eqn:linear relations} is 0 for all $R \in \CV_{-\bn}$. In the formulas above, $\tzeta_{ij} = \zeta_{ij}$ if $i \neq j$, but $\tzeta_{ii}$ denotes the Laurent polynomial in the numerator of \eqref{eqn:def zeta same}.

\end{proposition}

\medskip

\noindent It is easy to see that $\Phi \circ \Psi = 0$. Then Proposition \ref{prop:trivial} tells us that $R \in \CV_{-\bn}$ actually lies in $\CS_{-\bn}$ if and only if \eqref{eqn:linear relations} holds for a set of representatives
$$
\{ p_{\sigma}\}_{\sigma \in S(\bn)} \in \text{Ker }\Phi \Big/ \text{Im }\Psi
$$ 
The quotient $\text{Ker }\Phi / \text{Im }\Psi$ can be thought of as (a stronger version of) a first Koszul homology module, and it controls the linear relations that cut out $\CS_{-\bn} \subset \CV_{-\bn}$.

\medskip

\begin{proof} \emph{of Proposition \ref{prop:trivial}:} It suffices to fix $\sigma \neq \sigma'$. Then the claim of Proposition \ref{prop:trivial} boils down to the fact that
$$
\int_{|z_{\sigma(1)}| \gg \dots \gg |z_{\sigma(n)}|} f_{\sigma\sigma'}(z_1,\dots,z_n) R(z_1,\dots,z_n) \frac {\prod^{\sigma^{-1}(a) < \sigma^{-1}(b)}_{{\sigma'}^{-1}(a) > {\sigma'}^{-1}(b)} \left(1 - \frac {z_b}{z_a}\right)}{\prod^{\sigma^{-1}(a) < \sigma^{-1}(b)}_{{\sigma'}^{-1}(a) < {\sigma'}^{-1}(b)} \tzeta_{i_bi_a}\left(\frac {z_b}{z_a} \right)} \prod_{a=1}^n Dz_i = 
$$
$$
 = \int_{|z_{\sigma'(1)}| \gg \dots \gg |z_{\sigma'(n)}|} f_{\sigma\sigma'}(z_1,\dots,z_n) R(z_1,\dots,z_n) \frac {\prod^{\sigma^{-1}(a) < \sigma^{-1}(b)}_{{\sigma'}^{-1}(a) > {\sigma'}^{-1}(b)} \left(1 - \frac {z_b}{z_a}\right)}{\prod^{\sigma^{-1}(a) < \sigma^{-1}(b)}_{{\sigma'}^{-1}(a) < {\sigma'}^{-1}(b)} \tzeta_{i_bi_a}\left(\frac {z_b}{z_a} \right)} \prod_{a=1}^n Dz_i
$$
Indeed, the integrand is the same in both sides of the expression above, and there are no poles prohibiting us from moving the contours from $|z_a| \gg |z_b|$ if $\sigma^{-1}(a) > \sigma^{-1}(b)$ (as in the LHS) to $|z_a| \gg |z_b|$ if ${\sigma'}^{-1}(a) > {\sigma'}^{-1}(b)$ (as in the RHS).

\end{proof}

\medskip

\subsection{} 
\label{sub:words}

We will now set up the proof of Theorem \ref{thm:equal}. Consider the set of \textbf{letters}
$$
\left\{ i^{(d)} \right\}_{i \in I, d \in \BZ}
$$
A \textbf{word} is simply a sequence of letters
\begin{equation}
\label{eqn:word}
w = \left[ i_1^{(d_1)} \dots i_n^{(d_n)} \right]
\end{equation}
We will call $n$ the \textbf{length} of a word as above, and call
\begin{equation}
\label{eqn:exponents word}
\ow = (d_1,\dots,d_n)
\end{equation}
the sequence of exponents of $w$. The \textbf{degree} of the word \eqref{eqn:word} is defined as
\begin{equation}
\label{eqn:degree word}
\deg w = (\bs^{i_1} + \dots + \bs^{i_n}, d_1 + \dots + d_n) \in \nn \times \BZ
\end{equation}
Denote the set of all words by $\CW$. To $w \in \CW$ as above, we associate the elements
\begin{align}
&e_w = e_{i_1,d_1} \dots e_{i_n, d_n} \in \tUUp \label{eqn:associated word 1} \\
&f_w = e_{i_1,-d_1} \dots e_{i_n, -d_n} \in \tUUm \label{eqn:associated word 2}
\end{align}
It is clear that linear generating sets of $\oCS^\pm$ are given by
\begin{align}
&E_w := \tUpsilon^+(e_w) = \text{Sym} \left[ z_{i_1\bullet_1}^{d_1} \dots z_{i_n\bullet_n}^{d_n} \prod_{1\leq a < b \leq n} \zeta_{i_ai_b} \left(\frac {z_{i_a\bullet_a}}{z_{i_b\bullet_b}} \right) \right] \in \CV^+ \label{eqn:shuffle word 1} \\
&F_w := \tUpsilon^-(f_w) = \text{Sym} \left[ z_{i_1\bullet_1}^{-d_1} \dots z_{i_n\bullet_n}^{-d_n} \prod_{1\leq a < b \leq n} \zeta_{i_bi_a} \left(\frac {z_{i_b\bullet_b}}{z_{i_a\bullet_a}} \right) \right] \in \CV^- \label{eqn:shuffle word 2}
\end{align}
(in the formulas above, we choose $\{\bullet_a\}_{a\in \{1,\dots,n\}} \in \{1,\dots,n_{i_a}\}$ such that $\bullet_a < \bullet_b$ whenever $a < b$ and $i_a = i_b$) as $w$ runs over the set $\CW$ of all words.

\medskip

\begin{lemma} 
\label{lem:subset}

We have $\oCS^\pm \subseteq \CS^\pm$. 

\end{lemma}

\medskip

\begin{proof} We will prove the required statement for $\pm = -$, as the case of $\pm = +$ is analogous. The following formula will come in handy repeatedly
\begin{equation}
\label{eqn:pair e f}
\Big \langle e_v, F_w \Big \rangle = \Big \langle E_v, f_w \Big \rangle = 
\end{equation}
$$
= \int_{|z_1| \gg \dots \gg |z_n|} \mathop{\sum_{\sigma \in S(n)}}_{i_a  = j_{\sigma(a)}, \forall a} z_1^{d_1-k_{\sigma(1)}} \dots z_n^{d_n-k_{\sigma(n)}} \prod^{a<b}_{\sigma(a) > \sigma(b)} \frac {\zeta_{i_ai_b} \left(\frac {z_a}{z_b} \right)}{\zeta_{i_bi_a} \left(\frac {z_b}{z_a} \right)} \prod_{a=1}^n Dz_a
$$
for any pair of words
\begin{equation}
\label{eqn:two words}
v = \left [ i_1^{(d_1)} \dots i_n^{(d_n)} \right] \qquad \text{and} \qquad w = \left [ j_1^{(k_1)} \dots j_n^{(k_n)} \right]
\end{equation}
(the proof of \eqref{eqn:pair e f} is almost word-for-word as that of \cite[Remark 3.16]{Quiver 1}, so we will not repeat it here). Thus, the pairing $\langle \text{LHS of  \eqref{eqn:generator explicit}}, F_w \rangle$ equals
$$
\sum_{\sigma \in S(n)} \int_{|z_1| \gg \dots \gg |z_n|}  \frac {p_{j_{\sigma(1)},\dots,j_{\sigma(n)}}(z_1,\dots,z_n)}{z_1^{k_{\sigma(1)}} \dots z_n^{k_{\sigma(n)}}} \prod^{a<b}_{\sigma(a) > \sigma(b)} \frac {\zeta_{j_{\sigma(a)}j_{\sigma(b)}} \left(\frac {z_a}{z_b} \right)}{\zeta_{j_{\sigma(b)} j_{\sigma(a)}} \left(\frac {z_b}{z_a} \right)} \prod_{a=1}^n Dz_a
$$
In every summand over $\sigma$ in the formula above, let us change the variable according to $z_a = y_{\sigma(a)}$. Then the formula above reads
$$
\sum_{\sigma \in S(n)} \int_{|y_{\sigma(1)}| \gg \dots \gg |y_{\sigma(n)}|}  \frac {p_{j_{\sigma(1)},\dots,j_{\sigma(n)}}(y_{\sigma(1)},\dots,y_{\sigma(n)})}{y_1^{k_{1}} \dots y_n^{k_{n}}} \prod^{a<b}_{\sigma^{-1}(a) > \sigma^{-1}(b)} \frac {\zeta_{j_b j_a} \left(\frac {y_b}{y_a} \right)}{\zeta_{j_a j_b} \left(\frac {y_a}{y_b} \right)} \prod_{a=1}^n Dy_a
$$
As one moves the contours of integration from $|y_a| \gg |y_b| \Leftrightarrow \sigma^{-1}(a) < \sigma^{-1}(b)$ to $|y_a| \ll |y_b| \Leftrightarrow a<b$, one does not encounter any poles, so we conclude that
\begin{multline*}
\Big \langle \text{LHS of \eqref{eqn:generator explicit}}, F_w \Big \rangle = \\ 
= \int_{|y_1| \ll \dots \ll |y_n|} \sum_{\sigma \in S(n)}  \frac {p_{j_{\sigma(1)},\dots,j_{\sigma(n)}}(y_{\sigma(1)},\dots,y_{\sigma(n)})}{y_1^{k_{1}} \dots y_n^{k_{n}}} \prod^{a<b}_{\sigma^{-1}(a) > \sigma^{-1}(b)} \frac {\zeta_{j_b j_a} \left(\frac {y_b}{y_a} \right)}{\zeta_{j_a j_b} \left(\frac {y_a}{y_b} \right)} \prod_{a=1}^n Dy_a = 
\end{multline*}
$$
= \int_{|y_1| \ll \dots \ll |y_n|}  \frac {\sum_{\sigma \in S(n)}  p_{j_{\sigma(1)},\dots,j_{\sigma(n)}}(y_{\sigma(1)},\dots,y_{\sigma(n)}) \prod_{a<b} \zeta_{j_{\sigma(a)}j_{\sigma(b)}} \left(\frac {y_{\sigma(a)}}{y_{\sigma(b)}} \right)}{y_1^{k_{1}} \dots y_n^{k_{n}} \prod_{a<b}\zeta_{j_a j_b} \left(\frac {y_a}{y_b} \right)} \prod_{a=1}^n Dy_a 
$$
The numerator of the fraction above is precisely the LHS of \eqref{eqn:generator explicit equivalent}, which vanishes by assumption. We conclude that $\langle \text{LHS of \eqref{eqn:generator explicit}}, F_w \rangle = 0$. Since the LHS of \eqref{eqn:generator explicit} was an arbitrary element of $K^+$, we conclude that $F_w \in \CS^-$, as required.

\end{proof}

\begin{remark}
\label{rem:pairing}

Lemma \ref{lem:subset} implies that \eqref{eqn:descended pairing} and \eqref{eqn:descended pairing opposite} restrict to two pairings
\begin{equation}
\label{eqn:smallest pairing}
\oCS^+ \otimes \oCS^- \xrightarrow{\langle \cdot, \cdot \rangle} \BK
\end{equation}
which coincide due to the equality on the first line of \eqref{eqn:pair e f}. Moreover, because \eqref{eqn:descended pairing} and \eqref{eqn:descended pairing opposite} are non-degenerate in the $\CS^\pm$ argument, we conclude that \eqref{eqn:smallest pairing} is non-degenerate in both arguments. 

\end{remark}

\medskip

\subsection{} Let us fix a total order $<$ on the finite set $I$, and associate to it the following total order on the set of letters
$$
i^{(d)} < j^{(e)} \quad \text{if} \quad \begin{cases} d > e \\ \ \ \text{or} \\ d = e \text{ and } i<j \end{cases}
$$
Then we have the corresponding total lexicographic order on the set of words
$$
\left[i_1^{(d_1)} \dots i_n^{(d_n)} \right] < \left[j_1^{(e_1)} \dots j_m^{(e_m)} \right]
$$
if there exists $k$ such that $i_1^{(d_1)} = j_1^{(e_1)}, \dots, i_k^{(d_k)} = j_k^{(e_k)}$ and either $i_{k+1}^{(d_{k+1})} < j_{k+1}^{(e_{k+1})}$ or $k = n < m$.

\medskip

\begin{definition}
\label{def:non-increasing}

Consider $\{s_{ij} \in \BZ, \#_{ij} \geq 0 \}_{i,j\in I}$ uniquely defined so that
\begin{equation}
\label{eqn:explicit zeta ij}
\zeta_{ij}(x) = \alpha_{ij} \cdot x^{s_{ij}} + \dots + \beta_{ij} \cdot x^{s_{ij}+\#_{ij}}
\end{equation}
for all $i \neq j$, while
\begin{equation}
\label{eqn:explicit zeta ii}
\zeta_{ii}(x) = \frac {\alpha_{ii} \cdot x^{s_{ii}} + \dots + \beta_{ii} \cdot x^{s_{ii}+\#_{ii}}}{1-x}
\end{equation}
for all $i \in I$, where $\alpha_{ij}, \beta_{ij} \in \BK^\times$. Then a word \eqref{eqn:word} is called \textbf{non-increasing} if
\begin{equation}
\label{eqn:non-increasing}
\begin{cases} d_a < d_b - \sum_{a\leq x < b} (s_{i_x i_b} + s_{i_b i_x}) \\ \ \ \text{or} \\ d_a = d_b - \sum_{a\leq x < b} (s_{i_x i_b} + s_{i_b i_x}) \text{ and } i_a \geq i_b \end{cases}
\end{equation}
for all $1 \leq a < b \leq n$. Let $\CW_{\emph{non-inc}}$ denote the set of non-increasing words.

\end{definition}

\medskip

\begin{lemma}
\label{lem:finite}

There are finitely many non-increasing words of given degree, which are bounded above by any given word $v$. 

\end{lemma} 

\medskip

\begin{proof} Let us assume we are counting non-increasing words $[i_1^{(d_1)} \dots i_n^{(d_n)}]$ with $d_1+\dots+d_n = d$ for fixed $n$ and $d$. The fact that such words are bounded above implies that $d_1$ is bounded below. But then the inequality \eqref{eqn:non-increasing} implies that $d_2,\dots,d_n$ are also bounded below. The fact that $d_1+\dots+d_n$ is fixed implies that there can only be finitely many choices for the exponents $d_1,\dots,d_n$. Since there are also finitely many choices for $i_1,\dots,i_n \in I$, this concludes the proof.  

\end{proof}

\medskip

\subsection{}

Our main interest in the notion of non-increasing words is the following. 

\medskip

\begin{proposition} 
\label{prop:basis}

The set $\{e_w\}_{w \in \CW_{\emph{non-inc}}}$ is a $\BK$-linear basis of $\tUUp$.

\end{proposition}

\medskip

\begin{proof} Relation \eqref{eqn:rel quad} reads
\begin{multline*}
e_i(z) e_j(w) \left[ \alpha_{ji} \left(\frac wz \right)^{s_{ji}} + \dots + \beta_{ji} \left(\frac wz \right)^{s_{ji}+\#_{ji}} \right] = \\ = e_j(w) e_i(z) \left[ \alpha_{ij} \left(\frac zw \right)^{s_{ij}} + \dots + \beta_{ij} \left(\frac zw \right)^{s_{ij}+\#_{ij}} \right]
\end{multline*}
if $i \neq j$, and
\begin{multline*}
e_i(z) e_i(w) \left[\alpha_{ii} \left(\frac wz \right)^{s_{ii}} + \dots + \beta_{ii} \left(\frac wz \right)^{s_{ii}+\#_{ii}} \right] = \\ = e_i(w) e_i(z) \left[- \alpha_{ii} \left(\frac zw \right)^{s_{ii}-1} - \dots  - \beta_{ii} \left(\frac zw \right)^{s_{ii}+\#_{ii}-1} \right]
\end{multline*}
for all $i \in I$. Foiling out the formulas above and equating the coefficients of various monomials in $z$ and $w$ yields the equalities
$$
\alpha_{ji} \cdot e_{i,d} e_{j,k} + \text{larger} = \alpha_{ij} \cdot e_{j,k-s_{ij}-s_{ji}} e_{i,d+s_{ij}+s_{ji}} + \text{larger}
$$
if $i \neq j$, and
$$
\alpha_{ii} \cdot e_{i,d} e_{i,k} + \text{larger} = - \alpha_{ii} \cdot e_{i,k-2s_{ii}+1} e_{i,d+2s_{ii}-1} + \text{larger}
$$
where in either of the formulas above, the word ``larger" in the phrase ``$e_w+$larger" stands for linear combinations of $e_v$'s with $v > w$. Using the relations above, the proof of \cite[Proposition 3.11 in Subsection 4.4]{Quiver 1} explains how to write any $e_w$ for $w \notin \CW_{\text{non-inc}}$ as a linear combination of $e_v$'s with $v > w$ and $v \in  \CW_{\text{non-inc}}$. Following \loccitt, one actually obtains the stronger fact that there exists a constant $\beta(n)$ for all $n \in \BN$ (which only depends on the integers $\{s_{ij}, \#_{ij}\}_{i,j \in I}$) such that
\begin{equation}
\label{eqn:straighten}
e_w \in \sum^{v \in \CW_{\text{non-inc}} \text{ such that } v \geq w \text{ and}}_{\min \ow - \beta(\text{length}(w)) \leq \min \ov \leq \max \ov \leq \max \ow + \beta(\text{length}(w))} \BK \cdot e_v 
\end{equation}
for all $w \in \CW$. In other words, when writing any $e_w$ as a linear combination of $e_v$'s with $v \geq w$ and $v \in \CW_{\text{non-inc}}$, the exponents of any $v$ which appears in the linear combination will always remain within a global constant from the exponents of $w$.

\medskip

\noindent Let us now prove that $\{e_w\}_{w \in \CW_{\text{non-inc}}}$ are linearly independent in $\tUUp$. Given a monomial in variables $\{z_{i\bullet}\}_{i \in I, \bullet \in \BN}$, we will consider all ways to order its variables
\begin{equation}
\label{eqn:monomial}
\mu = z_{i_1\bullet_1}^{-k_1} \dots z_{i_n \bullet_n}^{-k_n}
\end{equation}
The \textbf{leading word} of $\mu$ as above is defined as the (lexicographically) largest word
\begin{equation}
\label{eqn:leading word}
w_\mu = \left [i_1^{(d_1)} \dots i_n^{(d_n)} \right] 
\end{equation}
among all ways to order the variables in \eqref{eqn:monomial}, where we write
\begin{equation}
\label{eqn:d and k}
d_a = k_a - \sum_{x > a} s_{i_x i_a} + \sum_{y < a}  s_{i_a i_y}, \qquad \forall a \in \{1,\dots,n\}
\end{equation}

\medskip

\begin{lemma} 
\label{lem:leading}

(\cite[Lemma 4.8]{Quiver 1}) The leading word of any monomial is non-increasing in the sense of Definition \ref{def:non-increasing}. 

\end{lemma}

\medskip

\noindent More generally, the leading word $\text{lead}(R)$ of any $R \in \CV^-$ is defined as the lexicographically largest of the leading words \eqref{eqn:leading word} for all the monomials which appear in $R$ with non-zero coefficient (the corresponding monomial will be called the \textbf{leading monomial} of $R$). Conversely, any non-increasing word $w$ is the leading word
\begin{equation}
\label{eqn:w as lead}
w = \text{lead}(\text{Sym } \mu)
\end{equation}
where the monomial $\mu$ is chosen such that $w_\mu = w$ in accordance with \eqref{eqn:leading word}. 

\medskip

\noindent Analogously to \cite[formula~(4.18)]{Quiver 1}, one can show by direct inspection that
\begin{equation}
\label{eqn:leading word pairing}
\Big \langle e_w, R \Big \rangle \text{ is } \begin{cases} \neq 0 &\text{if }w = \text{lead}(R) \\ = 0 &\text{if }w > \text{lead}(R) \end{cases}
\end{equation}
The formula above immediately shows the linear independence of the elements $e_w$, as $w$ runs over non-increasing words. Indeed, if one were able to write such an element $e_w$ as a linear combination of elements $\{e_v\}_{w < v \in \CW_{\text{non-inc}}}$, then this would contradict \eqref{eqn:leading word pairing} for $R = \text{Sym }\mu$ with $\mu$ as in \eqref{eqn:w as lead}. 

\end{proof}

\subsection{} We are now ready to conclude the proof of Theorem \ref{thm:equal}. \medskip

\begin{definition} 
\label{def:standard}

A word $w$ is called \textbf{standard} if 
\begin{equation}
\label{eqn:standard}
E_w \notin \sum_{v > w} \BK \cdot E_v
\end{equation}
Let $\CW_{\emph{stan}}$ denote the set of standard words. 
 
\medskip

\end{definition}

\noindent Because of the surjective algebra homomorphism $\widetilde{\Upsilon}^+ : \tUUp \rightarrow \oCS^+$, $e_w \mapsto E_w$, Proposition \ref{prop:basis} implies that any standard word is non-increasing. 

\medskip

\begin{proof} \emph{of Theorem \ref{thm:equal}:} We will prove the case of $\pm = -$, as that of $\pm = +$ is analogous. Fix any $n \in \BN$; all words in the present proof will have length $n$. As
\begin{equation}
\label{eqn:pole at infinity}
\frac {\zeta_{ij}(x)}{\zeta_{ji}(x^{-1})} \in \BK^\times \cdot x^{s_{ij}+s_{ji}+\#_{ij} - \delta_{ij}} + \left( \text{smaller deegree as }x \rightarrow \infty \right)
\end{equation}
for all $i,j \in I$, formula \eqref{eqn:pair e f} implies that $\langle e_v, F_w \rangle \neq 0$ (for words $v$ and $w$ as in \eqref{eqn:two words}) only if there exists $\sigma \in S(n)$ and natural numbers
$$
\Big\{ c_{ab} \leq s_{i_a i_b} + s_{i_b i_a} + \#_{i_a i_b} - \delta_{i_ai_b}\Big\}_{a < b, \sigma(a) > \sigma(b)}
$$
such that we have for all $a \in \{1,\dots,n\}$
\begin{equation}
\label{eqn:d and k graph}
d_a + \sum^{y > a}_{\sigma(y) < \sigma(a)} c_{ay} - \sum^{x<a}_{\sigma(x) > \sigma(a)} c_{xa} = k_{\sigma(a)}
\end{equation}
Thus, let us fix a natural number $m \gg$ the finitely many integers $\{s_{ij}, \#_{ij}\}_{i, j \in I}$ and $n$. Consider the infinite graph $G_{m,n}$ whose vertices are collections of integers
$$
(d_1,\dots,d_n) \text{ such that } d_a \leq d_{a+1}+m
$$
and whose edges are
\begin{equation}
\label{eqn:edge directed graph}
(d_1,\dots,d_n) \rightarrow ( k_1, \dots, k_n)
\end{equation}
if and only if there exists $\sigma \in S(n)$ and natural numbers
$$
\Big\{ c_{ab} \leq m \Big\}_{a < b, \sigma(a) > \sigma(b)}
$$
such that formula \eqref{eqn:d and k graph} holds (although the graph $G_{m,n}$ is defined to be directed, it is easy to see that an edge \eqref{eqn:edge directed graph} exists if and only if the opposite edge exists). The preceding discussion implies that
\begin{equation}
\label{eqn:only if}
\Big \langle e_v, F_w \Big \rangle \neq 0 \quad \Rightarrow \quad \text{exists edge } \ov \rightarrow \ow
\end{equation}
in $G_{m,n}$, for any two non-increasing words $v$ and $w$.

\medskip

\begin{lemma} (\cite[Lemma 3.18]{Quiver 1}) All connected components of $G_{m,n}$ are finite. \end{lemma}

\medskip

\noindent With the Lemma above in mind, let us define the following \underline{finite-dimensional} $\BK$-vector subspaces for any connected component $H \subset G_{m,n}$:
\begin{align*}
&\oCS^+_H = \sum_{w \in \CW_{\text{non-inc}}, \ \ow \in H} \BK \cdot E_w \subset \oCS^+ \\
&\oCS_H^- = \sum_{w \in \CW_{\text{non-inc}}, \ \ow \in H} \BK \cdot F_w \subset \oCS^-
\end{align*}
Because of \eqref{eqn:only if}, the restriction of the pairing \eqref{eqn:smallest pairing} satisfies
\begin{equation}
\label{eqn:pairing vanishes}
\Big \langle \oCS^+_H, \oCS_{H'}^- \Big \rangle = 0
\end{equation}
for any distinct connected components $H \neq H'$. The following is an immediate consequence of \eqref{eqn:pairing vanishes} and the non-degeneracy of the pairing \eqref{eqn:smallest pairing}.

\medskip

\begin{proposition}
\label{eqn:direct 1}

The restriction of the pairing \eqref{eqn:smallest pairing} to
\begin{equation}
\label{eqn:restricted pairing}
\oCS^+_H \otimes \oCS^-_{H} \xrightarrow{\langle \cdot, \cdot \rangle} \BK
\end{equation}
is non-degenerate for any connected component $H \subset G_{m,n}$.

\end{proposition}

\medskip

\begin{proposition}
\label{eqn:direct 2}

For any $n \in \BN$, we have 
\begin{equation}
\label{eqn:direct sum 1}
\bigoplus_{|\bn|=n} \oCS^\pm_{\bn} = \bigoplus^{\text{connected components}}_{H \subset G_{m,n}} \oCS^\pm_H
\end{equation}

\end{proposition}

\medskip

\begin{proof} Let us prove \eqref{eqn:direct sum 1} for $\pm = +$. Because the $E_w$'s span $\oCS^+$ as $w$ runs over all non-increasing words, all that we need to do to prove \eqref{eqn:direct sum 1} is to show that there are no linear relations among the various direct summands of the RHS. To this end, assume that we had a relation
$$
\sum_{H \text{ a connected component of } G_{m,n}} \alpha_H = 0
$$
for various $\alpha_H \in \oCS_H^+$. Pairing the relation above with a given $\oCS^-_H$ implies that
$$
\left \langle \alpha_H, \oCS^-_H \right \rangle = 0
$$
Because the pairing \eqref{eqn:restricted pairing} is non-degenerate, this implies that $\alpha_H = 0$. 

\end{proof}

\begin{proposition}
\label{prop:direct 3}

For any $n\in \BN$, we have
\begin{equation}
\label{eqn:direct sum 2}
\oCS^+_H = \bigoplus_{w \text{ standard}, \ \ow \in H} \BK \cdot E_w
\end{equation}
\begin{equation}
\label{eqn:direct sum 3}
\oCS^-_H = \bigoplus_{w \text{ standard}, \ \ow \in H} \BK \cdot F_w
\end{equation}
for any connected component $H \subset G_{m,n}$ \footnote{A priori, the notion of ``standard" words defined as in \eqref{eqn:standard} might differ from the analogous notion with $F$'s instead of $E$'s. The indexing sets defining $\oCS^+_H$ and $\oCS^-_H$ in \eqref{eqn:direct sum 2} and \eqref{eqn:direct sum 3} must therefore run over the notion of standard words defined with respect to $E$'s and $F$'s, respectively.}.

\end{proposition}

\medskip

\begin{proof} Since $\oCS^\pm_H$ are finite-dimensional vector spaces, the Proposition boils down to the general claim that any vector space spanned by vectors $\alpha_1,\dots,\alpha_k$ has a basis consisting of those $\alpha_i$'s which cannot be written as linear combinations of $\{\alpha_j\}_{j > i}$. 

\end{proof}

\medskip 

\noindent To complete the proof of Theorem \ref{thm:equal}, consider any $R \in \CS_{-\bn}$ with $|\bn| = n$. From \eqref{eqn:pair formula}, it is easy to see that:
$$
\Big \langle E_{[i_1^{(d_1)} \dots i_n^{(d_n)} ]}, R \Big \rangle = 0
$$
if $d_1$ is small enough. However, by Lemma \ref{lem:finite}, there are only finitely many non-increasing words $w$ of given degree with $d_1$ bounded below. This implies that
$$
\Big \langle E_w, R \Big \rangle \neq 0
$$
only for finitely many non-increasing words $w$. Let $H_1, \dots, H_t \subset G_{m,n}$ denote the connected components which contain the sequences of exponents of the aforementioned finitely many non-increasing words. The non-degeneracy of the pairings \eqref{eqn:restricted pairing} of finite-dimensional vector spaces implies that there exists an element
$$
R' \in \oCS^-_{H_1} \oplus \dots \oplus \oCS^-_{H_t} \subset \oCS^-
$$
such that $\langle E_w, R \rangle = \langle E_w, R' \rangle$ for all non-increasing words $w$. Then the non-degeneracy of \eqref{eqn:descended pairing} in the second argument implies that $R = R'$, as we needed to establish $\CS^- = \oCS^-$. The statement about the pairings \eqref{eqn:descended pairing} and \eqref{eqn:descended pairing opposite} being identical therefore follows from Remark \ref{rem:pairing}.

\end{proof}

\begin{remark}
\label{rem:construct}

Reversing the argument in the proof of Theorem \ref{thm:equal} allows us to find (algorithmically) for any $R \in \CV^- \backslash \CS^-$ some element $\phi \in K^+$ such that 
\begin{equation}
\label{eqn:construct}
\Big \langle \phi, R \Big \rangle \neq 0
\end{equation}
Indeed, one need only isolate the finitely many non-increasing words $w$ such that
$$
\Big \langle e_w, R \Big \rangle \neq 0
$$
and let $H_1,\dots,H_t \subset G_{m,n}$ denote the connected components which contain the sequences of exponents of the aforementioned finitely many non-increasing words. For any $w \in \CW_{\emph{non-inc}} \backslash \CW_{\emph{stan}}
$, we have by definition
$$
E_w = \sum_{v > w} \gamma_w^v \cdot E_v
$$
for certain $\gamma_w^v \in \BK$. Therefore,
$$
e_w - \sum_{v > w} \gamma_w^v \cdot e_v \in K^+
$$
and we claim that there exists $w \in \CW_{\emph{non-inc}} \backslash \CW_{\emph{stan}}$ such that
\begin{equation}
\label{eqn:not zero}
\Big \langle e_w - \sum_{v > w} \gamma_w^v \cdot e_v, R \Big \rangle \neq 0
\end{equation}
Indeed, if \eqref{eqn:not zero} failed for all $w \in \CW_{\emph{non-inc}} \backslash \CW_{\emph{stan}}$ then $R$ would pair trivially with the whole of $K^+$, and this would violate the assumption $R \in \CV^- \backslash \CS^-$. Moreover, one can achieve \eqref{eqn:not zero} for one of those finitely many $w$'s whose sequence of exponents $\ow$ lies in the connected components $H_1,\dots,H_t \subset G_{m,n}$, because all other $e_w$'s (as well as all $e_v$'s with $v > w$) pair trivially with $R$ by construction.

\end{remark}

\medskip

\subsection{} The quadratic quantum group of Definition \ref{def:quad} is merely the starting point of our study. The object we are actually interested in is the following.

\medskip

\begin{definition}
\label{def:quantum half}

The (positive/negative part of the) \textbf{quantum loop group} associated to the datum $\{\zeta_{ij}\}_{i,j \in I}$ is the $\BK$-algebra
\begin{equation}
\label{eqn:quantum half}
\UUpm = \tUUpm \Big/ K^\pm
\end{equation}
where $K^\pm$ denotes the kernel of the homomorphism $\tUpsilon^{\pm} : \tUUpm \longrightarrow \CV^\pm$.

\end{definition}

\medskip

\noindent As a direct consequence of the definition, the homomorphisms $\tUpsilon^\pm$ descend to isomorphisms
\begin{equation}
\label{eqn:upsilon}
\Upsilon^\pm : \UUpm \xrightarrow{\sim} \oCS^\pm = \CS^\pm
\end{equation}
with the last equality due to Theorem \ref{thm:equal}. The descended pairings \eqref{eqn:descended pairing}--\eqref{eqn:descended pairing opposite} therefore yield pairings
\begin{align}
&\UUp \otimes \CS^- \xrightarrow{\langle \cdot, \cdot \rangle} \BK  \label{eqn:quantum pairing 1} \\
&\CS^+ \otimes \UUm \xrightarrow{\langle \cdot, \cdot \rangle} \BK  \label{eqn:quantum pairing 2}
\end{align}
which both coincide under the isomorphisms \eqref{eqn:upsilon}
\begin{equation}
\label{eqn:quantum pairing}
\UUp \otimes \UUm \xrightarrow{\langle \cdot, \cdot \rangle} \BK 
\end{equation}
as shown in Remark \ref{rem:pairing}.

\medskip

\subsection{}
\label{sub:double}

In the present Subsection, let us assume that 
\begin{equation}
\label{eqn:symmetric}
\#_{ij} = -s_{ij} - s_{ji} + \delta_{ij}
\end{equation}
for all $i,j \in I$, in the notation of \eqref{eqn:explicit zeta ij}--\eqref{eqn:explicit zeta ii}. In particular, this implies $\#_{ij} = \#_{ji}$ for all $i,j$. We will refer to the collection $\{\zeta_{ij}\}_{i,j \in I}$ as being \textbf{symmetric}. Then all algebras studied in the present Section can be made into topological bialgebras, through the following procedure. First, we consider the extended algebras
\begin{align}
&\UUg = \frac {\UUp[h_{i,0}^{\pm 1}, h_{i,1}, h_{i,2}, \dots]_{i \in I}}{\left( h_i(z) e_j(w) = e_j(w) h_i(z) \frac {\zeta_{ij} \left(\frac zw \right)}{\zeta_{ji} \left(\frac wz \right)} \right)} \label{eqn:uug} \\
&\UUl = \frac {\UUm[{h'}_{i,0}^{\pm 1}, h'_{i,-1}, h'_{i,-2}, \dots]_{i \in I}}{\left( h'_i(z) f_j(w) = f_j(w) h'_i(z) \frac {\zeta_{ji} \left(\frac wz \right)}{\zeta_{ij} \left(\frac zw \right)} \right)} \label{eqn:uul}
\end{align}
where
$$
h_i(z) = \sum_{d=0}^{\infty} \frac {h_{i,d}}{z^d} \qquad \text{and} \qquad h'_i(z) = \sum_{d=0}^{\infty} h'_{i,-d} z^d
$$
Note that the assumption \eqref{eqn:symmetric} implies that
$$
\frac {\zeta_{ij} (x)}{\zeta_{ji} (x^{-1})}
$$
is regular and non-zero at $x = \infty$. With this in mind, one interprets the quotient relations in \eqref{eqn:uug} (respectively \eqref{eqn:uul}) by expanding them as power series in negative (respectively positive) powers of $\frac zw$. We can make $\UUg$ and $\UUl$ into topological bialgebras via the coproduct
\begin{align}
&\Delta(h_i(z)) = h_i(z) \otimes h_i(z) \label{eqn:coproduct h} \\
&\Delta(h'_i(z)) = h'_i(z) \otimes h'_i(z) \label{eqn:coproduct hh}
\end{align}
\begin{align}
&\Delta(e_i(z)) = h_i(z) \otimes e_i(z) + e_i(z) \otimes 1 \label{eqn:coproduct e} \\
&\Delta(f_i(z)) = 1 \otimes f_i(z) + f_i(z) \otimes h'_i(z) \label{eqn:coproduct f} 
\end{align}
It is straightforward to check that the pairing \eqref{eqn:quantum pairing} extends to a bialgebra pairing
\begin{equation}
\label{eqn:quantum pairing extended}
\UUg \otimes \UUl \xrightarrow{\langle \cdot, \cdot \rangle} \BK 
\end{equation}
via 
\begin{equation}
\label{eqn:pairing h}
\Big \langle h_i(z), h_i'(w) \Big \rangle = \frac {\zeta_{ij} \left(\frac zw \right)}{\zeta_{ji} \left(\frac wz \right)}  \qquad \text{expanded as } |z| \gg |w|
\end{equation}
and properties
\begin{align*}
&\Big \langle x, yy' \Big \rangle = \Big\langle \Delta(x), y \otimes y' \Big \rangle \\
&\Big \langle xx', y \Big \rangle = \Big\langle x' \otimes x, \Delta(y) \Big \rangle 
\end{align*}
for all $x,x' \in \UUg$ and $y,y' \in \UUl$. There are also antipode maps on $\UUg$ and $\UUl$ satisfying the usual properties in a topological Hopf algebra.

\medskip

\begin{definition}
\label{def:quantum}

The \textbf{quantum loop group} is defined as
\begin{equation}
\label{eqn:quantum}
\UU = \UUg \otimes \UUl \Big / \Big( h_{i,0} \otimes h_{i,0}' = 1 \otimes 1 \Big)
\end{equation}
with the multiplication governed by the Drinfeld double relation
\begin{equation}
\label{eqn:dd}
(x \otimes y) (x' \otimes y') = \sum_a \sum_b \Big \langle S(x'_{1,a}), y_{1,b} \Big \rangle x x'_{2,a} \otimes y_{2,b}  y'\Big \langle x'_{3,a}, y_{3,b} \Big \rangle
\end{equation}
for any $x,x' \in \UUg$ and $y,y' \in \UUl$ whose twice iterated coproducts satisfy
$$
\Delta^{(2)}(x') = \sum_a x'_{1,a} \otimes x'_{2,a} \otimes x'_{3,a} \qquad \text{and} \qquad \Delta^{(2)}(y) = \sum_b y_{1,b} \otimes y_{2,b} \otimes y_{3,b}
$$
The algebra $\UU$ inherits a Hopf algebra structure from its subalgebras $\UUg$ and $\UUl$.

\end{definition}

\medskip

\noindent Using relation \eqref{eqn:dd}, it is straightforward to deduce commutation relations between $e_i(z), h_i(z)$ and $f_j(w),h_j'(w)$. The most notable of these is the relation
\begin{equation}
\label{eqn:ef comm}
[e_i(z), f_j(w)] = \delta_{ij} \delta \left(\frac zw \right) \cdot \frac {h_i(z) - h_i'(w)}{q-q^{-1}}
\end{equation}
where $\delta(x) = \sum_{d \in \BZ} x^d$ is a formal series, and $\delta_{ij}$ is the Kronecker delta function. 

\medskip

\begin{remark} 

Choosing to have the scalar $q-q^{-1}$ in the denominator of \eqref{eqn:ef comm} is simply a matter of convention, to match the analogous constructions for usual quantum groups. However, to have \eqref{eqn:ef comm} indeed follow from \eqref{eqn:dd}, one needs to rescale the pairings \eqref{eqn:pair formula} and \eqref{eqn:pair formula opposite} by $(q^{-1}-q)^{-n}$. This modification does not substantially change any of the contents of the present Section, so we will ignore it.

\end{remark}

\bigskip

\section{The roots of the $\zeta$ functions}
\label{sec:roots}

\medskip

\noindent In the present Section, we will work over the splitting field of the Laurent polynomials $\zeta_{ij}(x) (1-x)^{\delta_{ij}}$ of \eqref{eqn:def zeta diff}--\eqref{eqn:def zeta same}, and determine the structure of the resulting shuffle algebras in terms of the pattern of zeroes of these Laurent polynomials.

\medskip

\subsection{}

In the present Section, we assume that 
\begin{equation}
\label{eqn:tzeta to zeta}
\zeta_{ij}(x) = \frac {\tzeta_{ij}(x)}{(1-x)^{\delta_{ij}}}
\end{equation}
with
\begin{equation}
\label{eqn:zeta factor}
\tzeta_{ij}(x) = \alpha_{ij}x^{s_{ij}} \prod_{e=1}^{\#_{ij}} (1-xq^{ij}_e)
\end{equation}
for various scalars $q_e^{ij} \in \BK^\times$ (in other words, we extend the ground field of Section \ref{sec:shuffle} to the splitting field of the Laurent polynomials that appear in \eqref{eqn:def zeta diff}--\eqref{eqn:def zeta same}). 

\medskip

\begin{example} 
\label{ex:k-ha}

When $\BK = \BQ$ and
$$
\zeta_{ij}(x) = (1-x)^{\#_{ij} - \delta_{ij}}
$$
for some matrix of non-negative integers $A = (\#_{ij})_{i,j \in I}$, we will use the notation
$$
\tUUpm_A, \UUpm_A, \CV^\pm_A, \CS^\pm_A
$$
for the $\BQ$-algebras defined in Section \ref{sec:shuffle}. Note that $\CV_A^+$ is the $K$-theoretic Hall algebra with 0 potential (\cite{KS, Pad 1}) of the quiver with vertex set $I$ and $\#_{ij}$ arrows from $i$ to $j$.

\end{example}

\medskip

\noindent We will soon show that the shuffle algebras from Example \ref{ex:k-ha} are the building blocks of shuffle algebras for general zeta functions \eqref{eqn:zeta factor}. Lemmas \ref{lem:basic 1} and \ref{lem:basic 2} hold in the generality of \eqref{eqn:def zeta diff}--\eqref{eqn:def zeta same}, but we will only use them in the setting of \eqref{eqn:zeta factor}.

\medskip

\begin{lemma}
\label{lem:basic 1}

For any $\bn = (n_i)_{i \in I} \in \nn$, we have
\begin{equation}
\label{eqn:shuffle ideal general}
\CS_{\pm \bn} = \frac {\CJ_{\bn}}{\Delta_{\bn}} \cap \CV_{\pm \bn}
\end{equation}
where $\CJ_{\bn} \subset \BK[z_{i1}^{\pm 1}, \dots, z_{in_i}^{\pm 1}]_{i \in I}$ denotes the ideal
\begin{equation}
\label{eqn:ideal general}
\CJ_{\bn} = \left( \prod_{(i,a) < (j,b)} \tzeta_{ij} \left(\frac {z_{ia}}{z_{jb}} \right) \right)_{\text{total order }<\text{ on } \{(i,a)| i \in I, a \in \{1,\dots,n_i\}\}}
\end{equation}
and
\begin{equation}
\label{eqn:delta general}
\Delta_{\bn} = \prod_{i \in I} \prod_{1 \leq a < b \leq n_i} (z_{ia} - z_{ib})
\end{equation}
In \eqref{eqn:shuffle ideal general} and henceforth, $\frac {\CJ_{\bn}}{\Delta_{\bn}}$ denotes the ideal quotient $(\CJ_{\bn} : (\Delta_{\bn}))$.

\end{lemma}

\medskip

\begin{proof} The inclusion $\subseteq$ in \eqref{eqn:shuffle ideal general} is an immediate consequence of $\CS_{\pm \bn} = \oCS_{\pm \bn}$, which was established in Theorem \ref{thm:equal}. To conclude the proof, we need to show that any 
$$
R \in \frac {\CJ_{\bn}}{\Delta_{\bn}} \cap \CV_{\pm \bn}
$$
actually lies in $\CS_{\pm \bn}$. The formula above implies that
\begin{equation}
\label{eqn:r}
R(z_{i1},\dots,z_{in_i})_{i \in I}  = \sum_{<} p_<(z_{i1},\dots,z_{in_i})_{i \in I}  \prod_{(i,a) < (j,b)} \zeta_{ij} \left(\frac {z_{ia}}{z_{jb}} \right) 
\end{equation}
for some Laurent polynomials $p_<$, where the sum goes over all total orders on the set $\{(i,a)| i \in I, a \in \{1,\dots,n_i\}\}$. However, the fact that $R \in \CV_{\pm \bn}$ implies that symmetrization (i.e. summing over all permutations of $z_{i1},\dots,z_{in_i}$, for each $i \in I$ separately) of $R$ has the effect of multiplying $R$ by a positive integer. As the symmetrization of the right-hand side of \eqref{eqn:r} clearly lies in $\oCS_{\pm \bn} = \CS_{\pm \bn}$, then so does $R$. 

\end{proof}

\medskip

\begin{example} 
\label{ex:wheel} 

In the setting of Example \ref{ex:k-ha}, Lemma \ref{lem:basic 1} implies that we can have $\CS_{A,\pm \bn} \neq \CV_{A,\pm \bn}$ only if there exist $i_1,\dots,i_k \in I$ such that $\#_{i_\bullet i_{\bullet+1}} > 0$ for all $\bullet \in \{1,\dots,k\}$ (let $i_{k+1}=i_1$), which is precisely the existence of a \textbf{wheel} in the language of Subsection \ref{sub:wheel}. Indeed, the absence of such a wheel would imply that $\#_{ii} = 0, \forall i$ and that there exists a total order $<$ on $I$ such that $\#_{ij} = 0$ if $i < j$. Then the total order on $\{(i,a)| i \in I, 1 \leq a \leq n_i\}$ which has $(i,a) < (j,b)$ if $i < j$ or if $i = j$ and $a<b$ would imply $1 \in \CJ_{\bn}$, and \eqref{eqn:shuffle ideal general} would give us $\CS_{A,\pm \bn} = \CV_{A,\pm \bn}$. 

\end{example}

\medskip

\begin{lemma}
\label{lem:basic 2}

Assume that one can partition the set
\begin{equation}
\label{eqn:partition}
I = I_1 \sqcup \dots \sqcup I_k
\end{equation}
such that $\#_{ij} = 0$ (i.e. $\tzeta_{ij}(x)$ is a unit in the ring $\BK[x^{\pm 1}]$) whenever $i \in I_s$ and $j \in I_t$ for $s \neq t$. Then for any $\bn \in \nn$, we have
\begin{equation}
\label{eqn:factor shuffle}
\CS_{\pm \bn} = \bigotimes_{t=1}^k \CS_{\pm \bn^t}^{(t)}
\end{equation}
as $\BK$-vector subspaces of
$$
\BK[z_{i1}^{\pm 1}, \dots, z_{in_i}^{\pm 1}]^{\esym}_{i \in I} = \bigotimes_{t=1}^k \BK[z_{i1}^{\pm 1}, \dots, z_{in_i}^{\pm 1}]^{\esym}_{i \in I_t}
$$
where $\CS^{\pm (t)}$ denotes the shuffle algebra defined with respect to $\{\zeta_{ij}\}_{i,j \in I_t}$, and $\bn^t \in \BN^{I_t}$ denotes the projection of $\bn \in \nn$ corresponding to the subset $I_t \subset I$.

\end{lemma}

\medskip

\begin{proof} By analogy with \eqref{eqn:ideal general}, for any $t \in \{1,\dots,k\}$ consider
$$
\CJ_{\bn^t}^{(t)} = \left( \prod_{(i,a) < (j,b)} \tzeta_{ij} \left(\frac {z_{ia}}{z_{jb}} \right) \right)_{\text{total order }<\text{ on } \{(i,a) | i \in I^t, a \in \{1,\dots,n_i\}\}}
$$
as an ideal in $\BK[z_{i1}^{\pm 1}, \dots, z_{in_i}^{\pm 1}]_{i \in I_t}$. The equality
\begin{equation}
\label{eqn:equality 1}
\CJ_{\bn} = \bigotimes_{t=1}^k \CJ_{\bn^t}^{(t)}
\end{equation}
is an easy consequence of the fact that every generator of the ideal in the left-hand side is a product of generators of the ideals in the right-hand side, and vice versa (this uses the fact that $\tzeta_{ij}$ is a unit if $i$ and $j$ lie in different parts of the partition \eqref{eqn:partition}). Since $\Delta_{\bn} = \prod_{t=1}^k \Delta_{\bn^t}$, it is an easy exercise to show that \eqref{eqn:equality 1} implies
\begin{equation}
\label{eqn:equality 2}
\frac {\CJ_{\bn}}{\Delta_{\bn}} = \bigotimes_{t=1}^k \frac {\CJ_{\bn^t}^{(t)}}{\Delta_{\bn^t}}
\end{equation}
(indeed, the $\supseteq$ inclusion is trivial, while $\subseteq$ follows from the fact that the polynomials $\Delta_{\bn^1},\dots,\Delta_{\bn^k}$ do not have any variables in common; see Claim \ref{claim:closed point} for a more complicated instance of this argument). Formulas \eqref{eqn:shuffle ideal general} and \eqref{eqn:equality 2} imply \eqref{eqn:factor shuffle}.

\end{proof}

\subsection{} 
\label{sub:point}

In the setting of the zeta functions \eqref{eqn:zeta factor}, consider any $\bn \in \nn$. For any
\begin{equation}
\label{eqn:point}
p = (p_{ia})_{i \in I, a \in \{1,\dots, n_i\}} \in (\BK^\times)^n \Big / \BK^\times
\end{equation}
(where $n = |\bn|$ and $\BK^\times$ acts on $(\BK^\times)^n$ by simultaneous rescaling) we define the set of equivalence classes
\begin{equation}
\label{eqn:equivalence}
I^p = \Big\{(i,a) \Big| i\in I, a \in \{1,\dots,n_i\}\Big\} \Big/  (i,a) \sim (j,b) \text{ if }i = j \text{ and }p_{ia} = p_{jb} 
\end{equation}
and the polynomial
\begin{equation}
\label{eqn:zeta factor partial}
\tzeta^p_{(i,a)(j,b)}(x) = \prod_{1 \leq e \leq \#_{ij} \text{ such that }q_e^{ij} = \frac {p_{jb}}{p_{ia}}}  (1-xq^{ij}_e)
\end{equation}
The definition above is designed to achieve two goals. Firstly
\begin{equation}
\label{eqn:divisibility}
\tzeta^p_{(i,a)(j,b)} \text{ divides } \tzeta_{ij}
\end{equation}
in $\BK[x^{\pm 1}]$. Secondly, define the matrix
$$
A^p = \Big\{ \#^p_{(i,a)(j,b)} \Big\}_{(i,a), (j,b) \in I^p} 
$$
where $\#^p_{(i,a)(j,b)}$ denotes the number of $e \in \{1,\dots,\#_{ij}\}$ such that
$$
q_e^{ij} = \frac {p_{jb}}{p_{ia}}
$$
(i.e. the number of linear factors in \eqref{eqn:zeta factor partial}). If we let $\bn^p \in \BN^{I^p}$ denote the vector whose $(i,a)$ entry is the cardinality of the equivalence class of $(i,a)$ in \eqref{eqn:equivalence}, then we have an isomorphism of vector spaces
\begin{equation}
\label{eqn:iso of shuffles}
\CS_{p, \pm \bn^p} \xrightarrow{\sim} \CS_{A^p, \pm \bn^p}, \qquad R( z_{ia} )  \mapsto R( z_{ia} p_{ia} )
\end{equation}
where $\CS^{\pm}_p$ is the shuffle algebra associated to $I^p$ and the functions \eqref{eqn:zeta factor partial}, and $\CS^\pm_A$ is the shuffle algebra defined in Example \ref{ex:k-ha}. The upshot is that, in degree $\bn^p$, the shuffle algebra defined with respect to the divisors \eqref{eqn:zeta factor partial} is simply a rescaling of the shuffle algebra in Example \ref{ex:k-ha} (which only depends on the choice of a quiver).

\medskip

\subsection{} The following result explains our earlier statement that all shuffle algebras are built out of those from Example \ref{ex:k-ha}.

\medskip

\begin{proposition}
\label{prop:intersection}

For any $\bn = (n_i)_{i \in I} \in \nn$, we have 
\begin{equation}
\label{eqn:intersection}
\CS_{\pm \bn} = \bigcap_{p \in (\BK^\times)^n  / \BK^\times} \CS_{p, \pm \bn^p} 
\end{equation}
where both sides are $\BK$-vector subspaces of $\BK[z_{i1}^{\pm 1}, \dots, z_{in_i}^{\pm 1}]_{i \in I}$. 

\end{proposition}

\medskip

\noindent Thus, formulas \eqref{eqn:iso of shuffles} and \eqref{eqn:intersection} reduce the study of the shuffle algebra $\CS_{\pm \bn}$ to the understanding of the shuffle algebras of Example \ref{ex:k-ha}. In fact, it is clear that the right-hand side of \eqref{eqn:intersection} is essentially a finite intersection, as $\CS_{p, \pm \bn^p}$ only depends on the equivalence relation defining \eqref{eqn:equivalence} and on the collection of linear factors that appear in \eqref{eqn:zeta factor partial}, both of which entail finitely many choices for each $\bn \in \nn$.

\medskip

\begin{proof} \emph{of Proposition \ref{prop:intersection}:} Let us recall the ideal \eqref{eqn:ideal general} and the polynomial \eqref{eqn:delta general}, and define similarly for any $p \in (\BK^\times)^n/\BK^\times$ the ideal
\begin{equation}
\label{eqn:ideal point}
\CJ_{p, \bn^p} = \left( \prod_{(i,a) < (j,b)} \tzeta_{(i,a)(j,b)}^p \left(\frac {z_{ia}}{z_{jb}} \right) \right)_{\text{total order }<\text{ on } \{(i,a)| i\in I, a \in \{1,\dots,n_i\}\}} 
\end{equation}
and the polynomial
\begin{equation}
\label{eqn:delta point}
\Delta_{\bn^p} = \prod_{i \in I} \mathop{\prod_{1 \leq a < b \leq n_i}}_{p_{ia} = p_{ib}} (z_{ia} - z_{ib}) 
\end{equation}
both in $\BK[z_{i1}^{\pm 1}, \dots, z_{in_i}^{\pm 1}]_{i \in I}$. Formula \eqref{eqn:shuffle ideal general} yields
\begin{equation}
\label{eqn:conclude 3}
\bigcap_{p \in (\BK^\times)^n  / \BK^\times} \CS_{p, \pm \bn^p}  = \bigcap_{p \in (\BK^\times)^n  / \BK^\times} \left( \frac {\CJ_{p,\bn^p}}{\Delta_{\bn^p}} \cap \BK[z_{i1}^{\pm 1}, \dots, z_{in_i}^{\pm 1}]^{\sym_p}_{i \in I} \right)
\end{equation}
where $\sym_p$ refers to polynomials which are symmetric in all $z_{ia}$ and $z_{ib}$ such that $p_{ia} = p_{ib}$. Thus, formula \eqref{eqn:intersection} is a consequence of
\begin{equation}
\label{eqn:remains to show}
\frac {\CJ_{\bn}}{\Delta_{\bn}} = \bigcap_{p \in (\BK^\times)^n  / \BK^\times} \frac {\CJ_{p,\bn^p}}{\Delta_{\bn^p}}
\end{equation}
which we will now prove.

\medskip

\begin{claim}
\label{claim:alg closed}

For any $p \in (\oBK^\times)^n  / \oBK^\times$, there exists $p' \in (\BK^{\times})^n / \BK^{\times}$ such that
\begin{equation}
\label{eqn:alg closed}
\CJ_{p',\bn^{p'}} = \CJ_{p,\bn^p} \qquad \text{and} \qquad \Delta_{\bn^{p'}} = \Delta_{\bn^p}
\end{equation}
where $\oBK$ denotes the algebraic closure of $\BK$. 

\end{claim}

\medskip

\noindent We will prove Claim \ref{claim:alg closed} by considering the partition
$$
\{(i,a)| i \in I, a \in \{1,\dots,n_i\}\} = \bigsqcup_{\text{cosets } H \text{ of }\oBK^\times/\BK^\times} \{(i,a) |  p_{ia} \in H\}
$$
For any coset $H$ of $\oBK^{\times} / \BK^{\times}$, let us simultaneously rescale $\{p_{ia} \in H\}$ to be generic elements of $\{p'_{ia} \in \BK^\times\}$. We claim that
\begin{equation}
\label{eqn:p and p'}
\tzeta^p_{(i,a)(j,b)}(x) = \tzeta^{p'}_{(i,a)(j,b)}(x)
\end{equation}
for all $(i,a)$ and $(j,b)$, which would imply that $p' \in (\BK^{\times})^n / \BK^{\times}$ satisfies \eqref{eqn:alg closed}. Indeed, the equality above is obvious for $(i,a)$ and $(j,b)$ lying in the same coset. Meanwhile, when $(i,a)$ and $(j,b)$ are in different cosets, both sides of \eqref{eqn:p and p'} are equal to 1: for the LHS, this is because $p_{jb}/p_{ia} \notin \BK^\times$, and thus cannot be equal to any of the roots $q_e^{ij} \in \BK^{\times}$. For the RHS, this is because the field $\BK$ of characteristic 0 is infinite, so the fact that $p'_{jb}$ and $p'_{ia}$ are generically chosen means that their ratio would not be equal to any one of the finitely many $q_e^{ij}$'s.

\medskip

\begin{claim}
\label{claim:closed point}

We have
\begin{equation}
\label{eqn:closed point}
\frac {\CJ_{\bn}}{\Delta_{\bn}} = \bigcap_{p \in (\oBK^\times)^n  / \oBK^\times} \frac {\CJ_{p,\bn^p}}{\Delta_{\bn^p}}
\end{equation}

\end{claim}

\medskip

\noindent It is clear that Claims \ref{claim:alg closed} and \ref{claim:closed point} imply \eqref{eqn:remains to show}, so it remains to prove the latter claim. To show that two ideals of $\BK[z_{i1}^{\pm 1}, \dots, z_{in_i}^{\pm 1}]_{i \in I}$ are equal, it suffices to show that they are equal in the localization at every closed point 
$$
p = (p_{ia})_{i \in I, a \in \{1,\dots,n_i\}} \in (\oBK^{\times})^n
$$
In the case at hand, this follows from \eqref{eqn:1} and \eqref{eqn:2} below. Firstly, we need
\begin{equation}
\label{eqn:1}
\left( \frac {\CJ_{\bn}}{\Delta_{\bn}} \right)_{p} = \left( \frac {\CJ_{p,\bn^p}}{\Delta_{\bn^p}} \right)_{p}
\end{equation}
which is a straightforward consequence of the fact that
$$
\tzeta_{ij} \left( \frac {z_{ia}}{z_{jb}} \right) = \tzeta^p_{(i,a)(j,b)} \left( \frac {z_{ia}}{z_{jb}} \right) \cdot \left( \text{unit in the localization at }p  \right)
$$
and
$$
\Delta_{\bn} = \Delta_{\bn^p} \cdot \left( \text{unit in the localization at }p \right)
$$
Secondly, we need to show that
\begin{equation}
\label{eqn:2}
\left( \frac {\CJ_{p,\bn^p}}{\Delta_{\bn^p}} \right)_{p} \subseteq \left( \frac {\CJ_{p',\bn^{p'}}}{\Delta_{\bn^{p'}}} \right)_{p}
\end{equation}
for any $p' \neq p$. To prove \eqref{eqn:2}, note the obvious fact that
$$
\tzeta^{p}_{(i,a)(j,b)} \left(\frac {z_{ia}}{z_{jb}} \right) \quad \text{is divisible by} \quad \tzeta^{p'}_{(i,a)(j,b)} \left(\frac {z_{ia}}{z_{jb}} \right)
$$
in the localization at $p$. Therefore, we have 
$$
\Big(\CJ_{p, \bn^p} \Big)_{p} \subseteq \Big(\CJ_{p', \bn^{p'}} \Big)_{p}
$$
Thus, to prove \eqref{eqn:2} it suffices to show that
\begin{equation}
\label{eqn:3}
\left( \frac {\CJ_{p', \bn^{p'}}}{\Delta_{\bn^{p}}} \right)_{p} = \left( \frac {\CJ_{p', \bn^{p'}}}{\Delta_{\bn^{p'}}} \right)_{p}
\end{equation}
Equality \eqref{eqn:3} can only fail due to linear factors $z_{ia} - z_{ib}$ which appear in $\Delta_{\bn^p}$ but do not appear in $\Delta_{\bn^{p'}}$, i.e. whenever $p_{ia} = p_{ib}$ but $p'_{ia} \neq p'_{ib}$. It suffices to show that for any such linear factor, we have
$$
\left( \frac {\CJ_{p', \bn^{p'}}}{z_{ia} - z_{ib}} \right)_p = \Big(\CJ_{p', \bn^{p'}} \Big)_{p}
$$
or in other words, that
\begin{equation}
\label{eqn:4}
(z_{ia} - z_{ib}) x \in \Big(\CJ_{p', \bn^{p'}} \Big)_{p} \quad \Rightarrow \quad x \in \Big (\CJ_{p', \bn^{p'}} \Big )_{p}
\end{equation}
Let us construct the graph $G$ with vertex set $\{(i',a')| i' \in I, a' \in \{1,\dots,n_{i'}\}\}$ and an edge between $(i',a')$ and $(j',b')$ if there exist linear factors in
$$
\tzeta^{p'}_{(i',a')(j',b')}\left(\frac {z_{i'a'}}{z_{j'b'}} \right)
$$
which are non-units in the localization at $p$. There can only be such an edge if
$$
\frac {p'_{i'a'}}{p'_{j'b'}} = \frac {p_{i'a'}}{p_{j'b'}}
$$
Just like in \eqref{eqn:equality 1}, we have
$$
\Big(\CJ_{p', \bn^{p'}} \Big)_{p} = \CJ_1 \cdots \CJ_k
$$
where each $\CJ_\bullet$ is an ideal whose generators only involve the variables $z_{ia}$ for $(i,a)$ lying in the $\bullet$-th connected component of $G$. The following result is an easy exercise, whose proof we leave to the reader.

\medskip

\begin{lemma}
\label{lem:product is intersection}

Assume that we have ideals $\CJ_1, \dots, \CJ_k$ in $\BK[w_1^{\pm 1}, \dots, w_n^{\pm 1}]$, such that we can partition the variables $w_1,\dots,w_n$ into disjoint sets $V_1, \dots, V_k$ with the generators of each $\CJ_\bullet$ being Laurent polynomials in the variables of $V_\bullet$. Then
$$
\CJ_1 \dots \CJ_k = \CJ_1 \cap \dots \cap \CJ_k
$$

\end{lemma}

\medskip

\noindent Lemma \ref{lem:product is intersection} applies equally well to ideals in the localization of $\BK[z_{i1}^{\pm 1}, \dots, z_{in_i}^{\pm 1}]_{i \in I}$ at $p$, so we have
$$
\Big(\CJ_{p', \bn^{p'}} \Big)_{p} = \CJ_1 \cap \dots \cap \CJ_k
$$
As $p_{ia} = p_{ib}$ and $p'_{ia} \neq p'_{ib}$, this implies that $(i,a)$ and $(i,b)$ are in different connected components of the graph $G$. Thus, to prove \eqref{eqn:4} it suffices to show that
$$
(z_{ia}-z_{ib}) x \in \CJ_\bullet \quad \Rightarrow \quad x \in \CJ_\bullet
$$
where $\CJ_{\bullet}$ is an ideal generated by polynomials that involve at most one of the variables $z_{ia}$ and $z_{ib}$. This is an obvious fact, whose proof we leave to the reader.

\end{proof}

\subsection{}

We will now consider the inclusion $\CS_{A,\pm \bn} \subset \CV_{A,\pm \bn}$ of Example \ref{ex:k-ha}, for any matrix $A$ with non-negative entries and any $\bn \in \nn$. Given elements $\bk, \bn \in \nn$, we will write $\b0 \leq \bk \leq \bn$ if $0 \leq k_i \leq n_i$ for all $i \in I$. In this case, if we have an ideal
$$
\CI_{\bk} \subseteq \BQ[z_{i1}^{\pm 1}, \dots, z_{ik_i}^{\pm 1}]_{i \in I}
$$
then we will write
\begin{equation}
\label{eqn:big int}
\CI_{\bk}^{(\bn)} = \bigcap_{C = (\text{distinct } c_1^{(i)}, \dots, c_{k_i}^{(i)} \in \{1,\dots, n_i\})_{i \in I}} \text{ideal generated by }\ph_C(\CI_{\bk})
\end{equation}
as an ideal in $\BQ[z_{i1}^{\pm 1}, \dots, z_{in_i}^{\pm 1}]_{i \in I}$, where $\ph_C$ denotes the ring homomorphism
\begin{equation}
\label{eqn:phi c}
\BQ[z_{i1}^{\pm 1}, \dots, z_{ik_i}^{\pm 1}]_{i \in I} \xrightarrow{\ph_C} \BQ[z_{i1}^{\pm 1}, \dots, z_{in_i}^{\pm 1}]_{i \in I}, \qquad z_{ia} \mapsto z_{ic_a^{(i)}}
\end{equation}

\medskip

\begin{proposition}
\label{prop:primary}

For any $A = \{\#_{ij} \geq 0\}_{i,j \in I}$, there exist homogeneous ideals
\begin{equation}
\label{eqn:ideal}
\CI_{A, \bn} \subseteq \BQ[z_{i1}^{\pm 1}, \dots, z_{in_i}^{\pm 1}]_{i \in I}
\end{equation}
with the quotient $\BQ[z_{i1}^{\pm 1}, \dots, z_{in_i}^{\pm 1}]_{i \in I} / \CI_{A,  \bn}$ supported on the small diagonal
\begin{equation}
\label{eqn:small diag}
\Big\{z_{ia} = z_{jb} \Big| \forall i,j \in I, a \in \{1,\dots,n_i\}, b \in \{1,\dots,n_j\} \Big\}
\end{equation}
and such that for all $\bn \in \nn$ we have
\begin{equation}
\label{eqn:primary}
\CS_{A,\pm \bn} = \CV_{A, \pm \bn} \bigcap_{\b0 \leq \bk \leq \bn} \CI_{A, \bk}^{(\bn)}
\end{equation}

\end{proposition}

\medskip

\noindent The choice of the ideals $\CI_{A, \bn}$ is not unique, although we will give a natural construction in \eqref{eqn:choice} below. The study of these ideals is warranted by the fact that Proposition \ref{prop:primary} reduces the study of shuffle algebras to the ideals $\CI_{A,  \bn}$. Moreover, as these ideals have finite codimension in any homogeneous degree \footnote{This is because the quotient $\BQ[z_{i1}^{\pm 1}, \dots, z_{in_i}^{\pm 1}]_{i \in I} / \CI_{A, \bn} $ is a graded coherent sheaf over a punctured affine line, and thus a finite-dimensional $\BQ$-vector space in any homogeneous degree.}, then the inclusion \eqref{eqn:ideal} is cut out by finitely many linear conditions in any homogeneous degree.

\medskip

\begin{proof} \emph{of Proposition \ref{prop:primary}:} Recall from \eqref{eqn:shuffle ideal general} that $\CS_{A, \pm \bn} = \frac {\CJ_{A, \bn}}{\Delta_{\bn}} \cap \CV_{A,\pm \bn}$, where
\begin{equation}
\label{eqn:shuffle ideal explicit} 
\CJ_{A, \bn} = \left( \prod_{(i,a) < (j,b)} (z_{ia} - z_{jb})^{\#_{ij}} \right)_{\text{total order }< \text{ on } \{(i,a)|i\in I, a\in \{1,\dots,n_i\}\}}
\end{equation}
and $\Delta_{\bn}$ is given by \eqref{eqn:delta general}. It suffices to construct the ideals \eqref{eqn:ideal} so that
\begin{equation}
\label{eqn:primary proof}
\frac {\CJ_{A, \bn}}{\Delta_{\bn}} = \bigcap_{\b0 \leq \bk \leq \bn} \CI_{A, \bk}^{(\bn)}
\end{equation}
and we will do this by induction on $\bn$. Let $\CI_{A,\b0} = (1)$ and assume that $\CI_{A, \bk}$ have been constructed for all $\bk < \bn$ (this is shorthand for $\b0 \leq \bk \leq \bn$ and $\bk \neq \bn$) and let us construct $\CI_{A, \bn}$. Formula \eqref{eqn:primary proof} for $\bn$ replaced by $\bk$ implies that
\begin{equation}
\label{eqn:inclusion 1}
\frac {\CJ_{A, \bk}}{\Delta_{\bk}} \subseteq \CI_{A, \bk} 
\end{equation}
However, comparing \eqref{eqn:shuffle ideal explicit} for $\bn$ and $\bk$ implies that
\begin{equation}
\label{eqn:inclusion 2}
\CJ_{A, \bn} \subseteq \CJ_{A, \bk}^{(\bn)}
\end{equation}
(see \eqref{eqn:big int} for the notation in the right-hand side), simply because the linear factors that generate the ideal $\CJ_{A, \bn}$ are divisible by the linear factors that generate the ideal $\CJ_{A, \bk}$ for any $\bk < \bn$. Thus, we have
\begin{equation}
\label{eqn:inclusion 3}
\frac {\CJ_{A, \bn}}{\Delta_{\bn}} \subseteq \frac {\CJ_{A, \bk}^{(\bn)}}{\Delta_{\bn}} = \left( \frac {\CJ_{A, \bk}}{\Delta_{\bk}} \right)^{(\bn)}
\end{equation}
where the equality is proved just like \eqref{eqn:3} (we leave the details as an exercise to the reader). Combining \eqref{eqn:inclusion 1} and \eqref{eqn:inclusion 3} yields
$$
\frac {\CJ_{A, \bn}}{\Delta_{\bn}} \subseteq \bigcap_{\b0 \leq \bk < \bn} \CI_{A, \bk}^{(\bn)}
$$

\begin{claim}
\label{claim:finite module}

The quotient
\begin{equation}
\label{eqn:finite module}
\left( \bigcap_{\b0 \leq \bk < \bn} \CI_{A, \bk}^{(\bn)} \right) \Big/ \left( \frac {\CJ_{A,\bn}}{\Delta_{\bn}} \right)
\end{equation}
is supported on the small diagonal $\{z_{ia} = z_{jb}|i,j \in I, a\in \{1,\dots,n_i\}, b \in \{1,\dots,n_j\}\}$. 

\end{claim}

\medskip

\noindent Let us first conclude the proof of Proposition \ref{prop:primary}, and then prove Claim \ref{claim:finite module}. Let
\begin{equation}
\label{eqn:m and n def}
\CM := \bigcap_{\b0 \leq \bk < \bn} \CI_{A, \bk}^{(\bn)} \supseteq \frac {\CJ_{A,\bn}}{\Delta_{\bn}} =: \CN
\end{equation}
and let $\fm$ denote the ideal $(z_{ia} - z_{jb})$ of the small diagonal $\text{Spec }\BQ[x^{\pm 1}]$. Consider the following descending sequence of ideals, as $k$ runs over $\BN$
\begin{equation}
\label{eqn:chain}
\CM \supseteq \dots \supseteq \CM \cap (\CN + \fm^k) \supseteq \dots \supseteq \CN
\end{equation}
By Claim \ref{claim:finite module}, the successive quotients of all the inclusions above are supported on the small diagonal. But since all ideals involved are homogeneous, the fiber of \eqref{eqn:chain} over the small diagonal is a chain of finite-dimensional $\BQ$-vector spaces tensored with $\BQ[x^{\pm 1}]$. By the Krull intersection theorem in the local ring at $\fm$, there exists $d$ large enough such that
\begin{equation}
\label{eqn:m and n}
\CN = \CM \cap (\CN + \fm^d)
\end{equation}
If we let 
\begin{equation}
\label{eqn:choice}
\CI_{A, \bn} = \frac {\CJ_{A,\bn}}{\Delta_{\bn}} + \fm^d
\end{equation}
for this large enough value of $d$, then \eqref{eqn:m and n} is precisely the required \eqref{eqn:primary proof}. 

\medskip

\begin{remark} Note that we have been quite profligate in our choice of the ideal \eqref{eqn:choice}, and other choices (such as those arising from the primary decomposition of the ideal $\frac {\CJ_{A, \bn}}{\Delta_{\bn}}$) might achieve \eqref{eqn:primary proof} for a larger ideal $\CI_{A, \bn}$ than the one of \eqref{eqn:choice}.

\end{remark}

\medskip

\noindent Let us prove Claim \ref{claim:finite module}. Consider any closed point $p = (p_{ia}) \in (\overline{\BQ}^{\times})^n$ with not all the $p_{ia}$'s equal to each other. Thus, the partition
$$
\{(i,a)| i \in I, a \in \{1,\dots,n_i\}\} = C_1 \sqcup \dots \sqcup C_k
$$
where 
$$
p_{ia} = p_{jb} \quad \Leftrightarrow \quad (i,a),(j,b) \in C_{\bullet} \text{ for some }\bullet \in \{1,\dots,k\}
$$ 
has the property that $k \geq 2$. By analogy with the proof of Claim \ref{claim:closed point}, we have 
$$
\left( \frac {\CJ_{A, \bn}}{\Delta_{\bn}} \right)_{p} = \prod_{\bullet=1}^{k} \ph_{C_{\bullet}} \left( \frac {\CJ_{A, \bn_{\bullet}}}{\Delta_{\bn_{\bullet}}} \right)_{p} 
$$
where $\bn_{\bullet}$ is the vector counting the number of variables in the part $C_{\bullet}$. Just like in Lemma \ref{lem:product is intersection}, one may show that
$$
\prod_{\bullet=1}^{k} \ph_{C_{\bullet}} \left( \frac {\CJ_{A, \bn_{\bullet}}}{\Delta_{\bn_{\bullet}}} \right)_{p}  = \bigcap_{\bullet=1}^{k} \ph_{C_{\bullet}} \left( \frac {\CJ_{A, \bn_{\bullet}}}{\Delta_{\bn_{\bullet}}} \right)_{p}
$$
Since \eqref{eqn:primary proof} is known to hold for all $\bn_{\bullet} < \bn$ by the induction hypothesis, we have
$$
\bigcap_{\bullet=1}^{k} \ph_{C_{\bullet}} \left( \frac {\CJ_{A, \bn_{\bullet}}}{\Delta_{\bn_{\bullet}}} \right)_{p} \supseteq \left( \bigcap_{\b0 \leq \bk < \bn} \CI_{A, \bk}^{(\bn)} \right)_{p}
$$
The three equations above prove that the quotient \eqref{eqn:finite module} is 0 in the localization at $p$. Since this holds for all $p$ outside of the small diagonal, we conclude Claim \ref{claim:finite module}.

\end{proof}

\medskip

\begin{proof} \emph{of Theorem \ref{thm:roots 1}:} immediately from \eqref{eqn:iso of shuffles}, \eqref{eqn:intersection} and \eqref{eqn:primary}. In more detail, we can define the ideal \eqref{eqn:ideal intro} inductively in $\bn$ by $\CI_{p,\b0} = (1)$, while for all $\bn > \b0$
\begin{equation}
\label{eqn:choice general}
\CI_{p, \bn} = \frac {\CJ_{ \bn}}{\Delta_{\bn}} + \fm_p^d
\end{equation}
where $\CJ_{\bn}$, $\Delta_{\bn}$ are defined in \eqref{eqn:ideal general}, \eqref{eqn:delta general} respectively, $\fm_p$ is the ideal of the small diagonal $(z_{ia}) \in p \BK^\times$ and the natural number $d$ is chosen large enough so that the following analogue of \eqref{eqn:m and n} holds
$$
\left( \frac {\CJ_{ \bn}}{\Delta_{\bn}} \right)_p = \left( \bigcap_{\b0 \leq \bk < \bn} \CI_{p, \bk}^{(\bn)} \right)_p \cap \left( \frac {\CJ_{\bn}}{\Delta_{\bn}} + \fm_p^d \right)_p
$$

\end{proof}

\subsection{} Let us now develop the dual (i.e. at the level of $\tUUpm$ instead of $\CV^\pm$, using \eqref{eqn:principle}) treatment of the preceding Subsections. The goal is to realize the inclusion
$$
\CS^\mp \subset \CV^\mp
$$
as pairing with certain elements of $K^\pm \subset \tUUpm$ under \eqref{eqn:pair} and \eqref{eqn:pair opposite}.

\medskip

\begin{proposition}
\label{prop:dual}

For any $A = \{\#_{ij} \geq 0\}_{i,j \in I}$, there exist finite sets
\begin{equation}
\label{eqn:dual}
W_{A,\pm \bn, \pm d} \subset K_A^\pm = \emph{Ker } \tUpsilon^\pm \subset \tUU_{A,\pm \bn, \pm d}
\end{equation}
for all $(\bn, d) \in \nn \times \BZ$, such that for all $R^\pm \in \CV^\pm_A$ we have
\begin{align}
&\Big \langle J_A^+ , R^- \Big \rangle = 0 \qquad \Leftrightarrow \qquad R^- \in \CS^-_A \label{eqn:dual 1} \\
&\Big \langle R^+ , J_A^- \Big \rangle = 0 \qquad \Leftrightarrow \qquad R^+ \in \CS_A^+ \label{eqn:dual 2}
\end{align}
where $J_A^\pm  = (W_{A,\pm \bn,\pm d})_{(\bn,d) \in \nn \times \BZ}$ as an ideal \footnote{In Subsection \ref{sub:j is k} we will prove that $J_A^\pm = K_A^\pm$.} of $\tUUpm_A$.

\end{proposition} 

\medskip

\begin{proof} We will deal with the case of $\pm = +$, as the case $\pm = -$ is analogous. We assume that \eqref{eqn:dual} have been constructed for all $\b0 \leq \bk < \bn$ such that \eqref{eqn:dual 1} and \eqref{eqn:dual 2} hold in all degrees $< \bn$, and let us perform the construction for $\bn$. Let
$$
J'_{A,\bn} = \bigoplus_{d \in \BZ} J'_{A,\bn,d}
$$
denote the degree $\bn$ component of the ideal $(W_{A, \bk, d})_{\b0 \leq \bk < \bn, d \in \BZ} \subset \tUUp_A$.

\medskip

\begin{claim}
\label{claim:ideal}

For any $R^- \in \CV_{A,-\bn}$, we have
\begin{equation}
\label{eqn:ideal claim}
\Big \langle J'_{A,\bn}, R^- \Big \rangle = 0 \qquad \Leftrightarrow \qquad R^- \in  \CV_{A,-\bn} \bigcap_{\b0 \leq \bk < \bn} \CI_{A, \bk}^{(\bn)}
\end{equation}
in the notation of \eqref{eqn:primary}.

\end{claim}

\medskip

\noindent Let us first show how Claim \ref{claim:ideal} allows us to complete the proof of Proposition \ref{prop:dual}. Let $\CM$ and $\CN$ be as in \eqref{eqn:m and n def}, and let us write $\CM_{d}$ and $\CN_{d}$ for their homogeneous degree $d$ components, for any $d \in \BZ$. The quotient $\CM_{d} / \CN_{d}$ is a finite dimensional $\BQ$-vector space. Formula \eqref{eqn:ideal claim} tells us that if $R^- \in \CV_{A,-\bn,-d}$, then
$$
\Big \langle J'_{A,\bn}, R^- \Big \rangle = 0 \qquad \qquad \Leftrightarrow \qquad \qquad R^- \in  \CV_{A,-\bn,-d} \cap \CM_{-d}
$$
Definition \ref{def:shuf} of $\CS^-$ reads
$$
\Big \langle K^+, R^- \Big \rangle = 0 \quad \ \ \Leftrightarrow \ \ \quad R^- \in  \CV_{A,-\bn,-d} \cap \CN_{-d} = \CS_{A,-\bn,-d}
$$
The finite-dimensionality of $\CM_{-d} / \CN_{-d}$ means that we can choose finitely many elements $\phi_1,\dots,\phi_N \in K^+$ such that
$$
\Big \langle J'_{A,\bn} + \sum_{\bullet = 1}^N \BQ \cdot \phi_\bullet, R^- \Big \rangle = 0 \qquad \ \Leftrightarrow \ \qquad R^- \in \CS_{A,-\bn,-d}
$$
In more detail, $\phi_1,\dots,\phi_N \in K^+$ are successively defined to pair non-trivially with finitely many elements of $\CV_{A,-\bn,-d} \cap (\CM_{-d} / \CN_{-d})$, so $\phi_1,\dots,\phi_N$ can be constructed algorithmically as in Remark \ref{rem:construct}. Letting $W_{A,\bn,d} = \{\phi_1,\dots,\phi_N\}$ yields \eqref{eqn:dual 1} in degree $\bn$, as required. 

\medskip

\begin{proof} \emph{of Claim \ref{claim:ideal}:} Fix any $\b0 \leq \bk < \bn$ and let us write $n = |\bn|$ and $k = |\bk|$. We will choose a relabeling $\bs^{i_1}+\dots+\bs^{i_n} = \bn$ such that $\bs^{i_{n-k+1}}+\dots+\bs^{i_n} = \bk$. Recall from \eqref{eqn:pair formula} that for any $R^- \in \CV_{-\bn}$, $\phi \in J_{A,\bk}$ and any $l_1,\dots, l_{n-k} \in \BZ$, the pairing
\begin{equation}
\label{eqn:pear}
\Big \langle e_{i_1,l_1} \dots e_{i_{n-k},l_{n-k}} \phi, R^-  \Big \rangle
\end{equation}
is computed as follows. Expand
\begin{equation}
\label{eqn:expand}
\frac {R^-(z_1,\dots,z_n)}{\prod_{\{1,\dots,n-k\} \ni a < b \in \{1,\dots,n\}} \zeta_{i_bi_a} \left(\frac {z_b}{z_a} \right)}
\end{equation}
as a power series in $|z_1| \gg \dots \gg |z_{n-k}|$, extract the coefficient of $z_1^{-l_1}\dots z_{n-k}^{-l_{n-k}}$, and then pair the resulting Laurent polynomial in $z_{n-k+1},\dots,z_n$ with $\phi$. Since $l_1,\dots, l_{n-k}$ can be arbitrary, formula \eqref{eqn:dual 1} in degree $\bk$ means that the vanishing of all the pairings \eqref{eqn:pear} is equivalent to $R^-$ lying in the ideal generated by
$$
\ph_C(\CI_{A,\bk})
$$
where $\ph_C : \BQ[z_1^{\pm 1}, \dots z_k^{\pm 1}] \longrightarrow \BQ[z_1^{\pm 1}, \dots z_n^{\pm 1}]$ is the map that sends $z_i \mapsto z_{n-k+i}$. Since $R^-$ is symmetric, this is equivalent to $R^-$ lying in the intersection 
$$
\CI_{A,\bk}^{(\bn)}
$$
in the notation of \eqref{eqn:big int}. Finally, because $\b0 \leq \bk < \bn$ is arbitrary, then we obtain precisely the equivalence \eqref{eqn:ideal claim}.

\end{proof}

\end{proof}

\begin{remark}
\label{rem:shift} 

The proof given above is invariant under the shift homomorphisms \eqref{eqn:shift 1}, \eqref{eqn:shift 2}, \eqref{eqn:shift 3}, \eqref{eqn:shift 4}, in the sense that for fixed $\bk \in \zz$ one can choose
\begin{equation}
\label{eqn:shift w}
W_{A,\pm \bn, \pm (d + \sum_{i \in I} n_i k_i)} = \tau_{\bk} (W_{A,\pm \bn, \pm d})
\end{equation}
for any $(\bn, d) \in \nn \times \BZ$. Thus, if $k_i > 0$ for all $i \in I$, then equality \eqref{eqn:shift w} reduces the choice of $W_{A,\pm \bn, \pm d}$ to finitely many choices for each $\bn \in \nn$.  

\end{remark}

\medskip

\begin{remark}
\label{rem:ideal} 

We have been intentionally vague concerning which kind of ideal one considers in \eqref{eqn:dual 1}--\eqref{eqn:dual 2}. The proof given above works if $J^\pm_A$ is defined as the left ideal generated by the elements $W_{A,\pm \bn, \pm d}$, but a similar argument works if $J^\pm_A$ is defined as the right (or two-sided) ideal generated by the same elements.

\end{remark}

\medskip

\subsection{} Just like the study of shuffle algebras can be reduced to the (purely combinatorial) particular case of Example \ref{ex:k-ha}, the same can be said about the dual situation of quantum loop groups. Specifically, in the following Proposition we will show how to define the elements \eqref{eqn:elements intro} for the quantum loop group associated to any choice of zeta functions \eqref{eqn:zeta factor}, in terms of the finite sets 
$$
W_{A, \pm \bn, \pm d} \subset \tUU_{A, \pm \bn, \pm d}
$$
that we defined in \eqref{eqn:dual}, for various matrices $A = \{\#_{ij} \geq 0\}_{i,j \in I}$.

\medskip

\begin{definition}
\label{def:reduce}

For any $(\bn,d) \in \nn \times \BZ$ and $p \in (\BK^\times)^n/\BK^\times$, assume that
$$
W_{A^p,\bn^p, d} = \Big\{\phi_1,\dots,\phi_N \Big\} \subset \tUU_{A^p, \bn^p, d}
$$
where $A^p$ and $\bn^p$ are defined as in Subsection \ref{sub:point}. Let us choose a representation
$$
\phi_\bullet = \mathop{\sum_{\{(i_1,a_1), \dots, (i_n,a_n)\} =}}_{\{(i,a)| i\in I, a \in \{1,\dots,n_i\}\}} \Big [ \rho^\bullet_{(i_1,a_1),\dots,(i_n,a_n)} (z_1,\dots,z_n) e_{(i_1,a_1)}(z_1) \dots e_{(i_n,a_n)}(z_n) \Big]_{\emph{ct}} 
$$
for each $\bullet \in \{1,\dots,N\}$, where $\rho^\bullet_{(i_1,a_1),\dots,(i_n,a_n)}$ are Laurent polynomials. Define
$$
W_{p, \bn, d} = \Big \{\phi_1', \dots, \phi'_N \Big \}  \subset \tUU_{ \bn, d}
$$
given by the following formula for each $\bullet \in \{1,\dots,N\}$
\begin{multline}
\label{eqn:elements}
\phi'_\bullet = \mathop{\sum_{\{(i_1,a_1), \dots, (i_n,a_n)\} =}}_{\{(i,a)| i\in I, a \in \{1,\dots,n_i\}\}} \Big [ \rho^\bullet_{(i_1,a_1),\dots,(i_n,a_n)} (z_1,\dots,z_n) \\ \prod_{1\leq u < v \leq n} \frac {\tzeta_{i_vi_u} \left(\frac {z_v}{z_u} \right) \left( - \frac {z_u}{z_v} \right)^{\delta_{i_ui_v} \delta_{v \lhd u}}}{\tzeta^p_{(i_v,a_v)(i_u,a_u)} \left(\frac {z_v}{z_u} \right)} \cdot e_{i_1}(z_1) \dots e_{i_n}(z_n) \Big]_{\emph{ct}} 
\end{multline}
In the right-hand side of \eqref{eqn:elements}, we fix an arbitrary total order $\lhd$ on $I^p$ and write
$$
u \lhd v \qquad \text{if} \qquad (i_u,a_u) \lhd (i_v,a_v) 
$$
The analogous construction with $e \leftrightarrow f$ and $\tzeta_{i_vi_u}\left(\frac {z_v}{z_u} \right) \leftrightarrow \tzeta_{i_ui_v}\left(\frac {z_u}{z_v} \right)$ defines $W_{p,-\bn,-d}$.

\end{definition}

\medskip

\begin{proposition}
\label{prop:reduce}

For any $(\bn,d) \in \nn \times \BZ$ and $p \in (\BK^\times)^n/\BK^\times$, we have $W_{p, \pm \bn, \pm d} \in K^\pm = \emph{Ker }\tUpsilon^\pm$. Moreover, for all $R^\pm \in \CV^\pm$ we have
\begin{align}
&\Big \langle J^+ , R^- \Big \rangle = 0 \qquad \Leftrightarrow \qquad R^- \in \CS^- \label{eqn:dual 3} \\
&\Big \langle R^+ , J^- \Big \rangle = 0 \qquad \Leftrightarrow \qquad R^+ \in \CS^+ \label{eqn:dual 4}
\end{align}
where $J^\pm = (W_{p,\pm \bn, \pm d})_{p \in (K^\times)^n /\BK^\times, \bn \in \nn, d \in \BZ}$ as an ideal \footnote{In Subsection \ref{sub:j is k} we will prove that $J^\pm = K^\pm$.} of $\tUUp$.

\end{proposition}

\medskip

\begin{proof} We will deal with the case $\pm = +$, as the case $\pm = -$ is analogous and left as an exercise to the reader. Let us show that the elements \eqref{eqn:elements} lie in $K^+ = \text{Ker }\tUpsilon^+$. To keep the subsequent formulas short, we will suppress the indexing sets of $\sum$ and $\rho$. As we saw in \eqref{eqn:generator explicit} and \eqref{eqn:generator explicit equivalent}, $\tUpsilon^+$ applied to the right-hand side of \eqref{eqn:elements} is
$$
\Sym \left[ \sum \rho (z_1,\dots,z_n)  \prod_{1\leq u < v \leq n} \frac {\tzeta_{i_vi_u} \left(\frac {z_v}{z_u} \right) \left( - \frac {z_u}{z_v} \right)^{\delta_{i_ui_v} \delta_{v \lhd u}} \zeta_{i_ui_v} \left(\frac {z_u}{z_v} \right)}{\tzeta^p_{(i_v,a_v)(i_u,a_u)} \left(\frac {z_v}{z_u} \right)} \right]
$$
where ``Sym" denotes symmetrization with respect to all $z_u$ and $z_v$ such that $i_u = i_v$ (and one needs to perform a relabeling \eqref{eqn:relabeling} in order for the Laurent polynomial above to be an element of $\CV^+$). The expression above is equal to the symmetric rational function
$$
\prod_{1\leq u \neq v \leq n} \zeta_{i_ui_v} \left(\frac {z_u}{z_v} \right)
$$
times 
$$
\Sym \left[ \sum \rho(z_1,\dots,z_n)  \prod_{1\leq u < v \leq n} \frac {\left(1 - \frac {z_v}{z_u} \right)^{\delta_{i_u i_v}\delta_{(i_u,a_u) \not \sim (i_v,a_v)}}\left( - \frac {z_u}{z_v} \right)^{\delta_{i_ui_v} \delta_{v \lhd u}}}{\zeta^p_{(i_v,a_v)(i_u,a_u)} \left(\frac {z_v}{z_u} \right)} \right]
$$
where $\sim$ denotes the equivalence relation that defines $I^p$, see \eqref{eqn:equivalence}, and $\delta_{x \not \sim y}$ is equal to 1 if $x \not \sim y$ and 0 otherwise. The expression directly above is equal to 
$$
\mathop{\prod_{1 \leq u \neq v \leq n}}_{(i_u,a_u) \lhd (i_v,a_v)} \left(1 - \frac {z_v}{z_u} \right)^{\delta_{i_u i_v}}
$$
(which takes the same value for all the summands in $\sum$, after relabeling) times 
$$
\text{Sym} \left[ \sum \rho (z_1,\dots,z_n)  \prod_{1\leq u < v \leq n} \frac 1{\zeta^p_{(i_v,a_v)(i_u,a_u)} \left(\frac {z_v}{z_u} \right)} \right]
$$
The symmetrization above is 0 because $\phi_\bullet \in K_{A^p}^+$. By the chain of equalities above, we conclude that $\phi_\bullet' \in K^+$, as required.

\medskip

\noindent We have just proved that $J^+ \subseteq K^+$, and have therefore established the $\Leftarrow$ implication of \eqref{eqn:dual 3}. As for the $\Rightarrow$ implication, let us consider an element $R^- \in \CV^-$ such that $\langle J^+, R^- \rangle = 0$. First of all, note that for $\phi'_\bullet$ as in \eqref{eqn:elements}, we have
\begin{multline*}
\Big \langle \phi_\bullet', R^- \Big \rangle = \sum \int_{|z_1| \gg \dots \gg |z_n|} \rho (z_1,\dots,z_n) \\ \prod_{1\leq u < v \leq n} \frac {\tzeta_{i_vi_u} \left(\frac {z_v}{z_u} \right) \left( - \frac {z_u}{z_v} \right)^{\delta_{i_ui_v} \delta_{v \lhd u}}}{\tzeta^p_{(i_v,a_v)(i_u,a_u)} \left(\frac {z_v}{z_u} \right)} \cdot \frac {R^-(z_1,\dots,z_n)}{\prod_{1 \leq u < v \leq n} \zeta_{i_vi_u} \left(\frac {z_v}{z_u} \right)}
\end{multline*}
\begin{align*}
&=\sum \int_{|z_1| \gg \dots \gg |z_n|} \rho (z_1,\dots,z_n)  \frac {R^-(z_1,\dots,z_n)\prod_{1 \leq u \neq v \leq n} \left(1 - \frac {z_v}{z_u} \right)^{\delta_{i_u i_v}\delta_{(i_u,a_u) \lhd (i_v,a_v)}}
}{\prod_{1 \leq u < v \leq n} \zeta^p_{(i_v,a_v)(i_u,a_u)} \left(\frac {z_v}{z_u} \right)} \\
&= \left \langle \phi_\bullet, R^-(z_1,\dots,z_n) \mathop{\prod_{1 \leq u \neq v \leq n}}_{(i_u,a_u) \lhd (i_v,a_v)} \left(1 - \frac {z_v}{z_u} \right)^{\delta_{i_u i_v}}
\right \rangle
\end{align*}
where the pairing on the top line denotes the one between $\tUUp$ and $\CV^-$, while the pairing on the bottom line denotes the one between $\tUU^+_{A^p}$ and $\CV^-_{A^p}$. Analogous formulas as above hold when $\phi_\bullet$ and $\phi_\bullet'$ are multiplied by various products of $e_{i,d}$'s, thus yielding general elements of $J^+_{A^p}$ and $J^+$, respectively. The fact that $R^-$ pairs trivially with the whole of $J^+$ implies that
\begin{equation}
\label{eqn:tilde r}
\widetilde{R}^- = R^- \mathop{\prod_{1 \leq u \neq v \leq n}}_{(i_u,a_u) \lhd (i_v,a_v)} \left(1 - \frac {z_v}{z_u} \right)^{\delta_{i_u i_v}}
\end{equation}
pairs trivially with $J^+_{A^p}$, which by \eqref{eqn:dual 1} means that $\widetilde{R}^- \in \CS_p^-$. By \eqref{eqn:primary}, $\CS_p^-$ is the intersection of ideals $\CI$ which all enjoy the following analogue of property \eqref{eqn:4}
$$
(z_u - z_v) x \in \CI \quad \Rightarrow \quad x \in \CI
$$
whenever $i_u = i_v$ but $(i_u,a_u) \not \sim (i_v,a_v) \Leftrightarrow p_{i_ua_u} \neq p_{i_va_v}$ (indeed, the only non-trivial situation is when $\CI$ is generated by polynomials which involve both $z_u$ and $z_v$; but in this case, $\CI$ is supported on the locus 
$$
L = \left\{\frac {z_u}{z_v} = \frac {p_{i_ua_u}}{p_{i_va_v}} \neq 1 \right\}
$$
and $z_u - z_v$ is nowhere vanishing on $L$). Therefore, we conclude that $R^- \in \CS_p^-$; since this holds for all $p \in (\BK^\times)^n / \BK^\times$, \eqref{eqn:intersection} implies that $R^- \in \CS^-$, as required. 

\end{proof}

\medskip

\subsection{} 
\label{sub:j is k}

Now that we have defined the elements $W_{p,\pm \bn, \pm d} \in \tUUpm$ for any point \eqref{eqn:point} and any $(\bn,d) \in \nn \times \BZ$, it remains to show that the ideal $J^\pm$ they generate is equal to $K^\pm$, and this will imply Theorem \ref{thm:roots 2}. Comparing \eqref{eqn:shuf 2} with \eqref{eqn:dual 3} shows that
\begin{equation}
\label{eqn:not quite}
\Big \langle K^+, R^- \Big \rangle = 0 \quad \Leftrightarrow \quad \Big \langle J^+, R^- \Big \rangle = 0
\end{equation}
If $\tUUp \otimes \CV^- \xrightarrow{\langle \cdot, \cdot \rangle} \BK$ were a non-degenerate pairing of finite-dimensional vector spaces (or at least finite-dimensional in any $\zz \times \BZ$ degree), then \eqref{eqn:not quite} would imply that $J^+ = K^+$ and the proof of Theorem \ref{thm:roots 2} would be complete. Thus, we need to adapt this line of reasoning to the setting of infinite-dimensional vector spaces. 

\medskip

\begin{proof} \emph{of Theorem \ref{thm:roots 2}:} We will prove that $J^+ = K^+$, as the situation when $+$ is replaced by $-$ is analogous. Fix $(\bn, d) \in \nn \times \BZ$ and our goal will be to show that
\begin{equation}
\label{eqn:goal}
J_{\bn, d} = J^+ \cap \tUU_{\bn, d} \qquad \text{equals} \qquad K_{\bn,d} = K^+ \cap \tUU_{\bn,d}
\end{equation}
By construction, we have
\begin{equation}
\label{eqn:goal 1}
J_{\bn,d} \subseteq K_{\bn,d}
\end{equation}
and it remains to prove the opposite inclusion. For any finite set $T \subset \CW_{\text{non-inc}}$ of non-increasing words of degree $(\bn,d)$, define
$$
\tUU^{+,T} = \bigoplus_{w \in T} \BK \cdot e_w
$$
(recall that as $w$ runs over $\CW_{\text{non-inc}}$, the $e_w$'s yield a basis of $\tUUp$, as shown in Proposition \ref{prop:basis}). Moreover, the notion of ``leading word" introduced in the latter part of the proof of Proposition \ref{prop:basis} allows us to define
\begin{equation}
\label{eqn:v finite}
\CV^{-,T} \subset \CV_{-\bn,-d}
\end{equation}
to be the subspace of Laurent polynomials $R^-$, all of whose monomials have leading word in $T$. Thus, the restriction of the pairing \eqref{eqn:pair}
\begin{equation}
\label{eqn:finite pairing}
\tUU^{+,T} \otimes \CV^{-,T} \xrightarrow{\langle \cdot, \cdot \rangle} \BK
\end{equation}
is a non-degenerate pairing of finite-dimensional vector spaces (indeed, because the two vector spaces have the same dimension, it suffices to show non-degeneracy in the second argument, which follows immediately from \eqref{eqn:leading word pairing}). 

\medskip

\noindent Let us consider a finite set of non-increasing words $T$ of degree $(\bn,d)$, which will be chosen in \eqref{eqn:choice of t}. The only thing we postulate for the time being is that the set $T$ can be chosen ``arbitrarily large", i.e. to contain any given finite set of words. Let $\CS^{-,T} = \CS^- \cap \CV^{-,T}$. Formula \eqref{eqn:dual 3} implies that
$$
\CS^{-,T} \subseteq \left(J^+ \cap \tUU^{+,T} \right)^\perp
$$
where the orthogonal complement is defined with respect to the pairing \eqref{eqn:finite pairing}. Our goal will be to prove the opposite inclusion, namely the following result.

\medskip

\begin{claim}
\label{claim:claim}

We have:
\begin{equation}
\label{eqn:opposite inclusion}
\left(J^+ \cap \tUU^{+,T}\right)^\perp \subseteq \CS^{-,T}
\end{equation}

\end{claim}

\medskip

\noindent Let us first use \eqref{eqn:opposite inclusion} to establish the opposite inclusion to \eqref{eqn:goal 1}, and thus conclude the proof of Theorem \ref{thm:roots 2}. Because \eqref{eqn:finite pairing} is a non-degenerate pairing of finite-dimensional vector spaces, we have
$$
\CS^{-,T} = \left(J^+ \cap \tUU^{+,T}\right)^\perp  \qquad \Rightarrow \qquad J^+ \cap \tUU^{+,T} = \left(\CS^{-,T}\right)^\perp
$$
However, we have $\left(\CS^{-,T}\right)^\perp \supseteq K^+ \cap \tUU^{+,T}$ by \eqref{eqn:shuf 2}, so the display above implies
$$
J^+ \cap \tUU^{+,T} \supseteq K^+ \cap \tUU^{+,T}
$$
As the finite set $T$ may be chosen arbitrarily large, the formula above implies
\begin{equation}
\label{eqn:goal 2}
J_{\bn, d} \supseteq K_{\bn,d}
\end{equation}
for all $(\bn,d) \in \nn \times \BZ$, which precisely provides the opposite inclusion to \eqref{eqn:goal 1}. 

\medskip

\begin{proof} \emph{of Claim \ref{claim:claim}:} The statement we need to prove is that for any $R^- \in \CV^{-,T}$
\begin{equation}
\label{eqn:hard}
\Big \langle \phi, R^- \Big \rangle = 0, \ \forall \phi \in J^+ \cap \tUU^{+,T} \qquad \Rightarrow \qquad R^- \in \CS^{-,T}
\end{equation}
Equivalently, we will prove its contrapositive statement, which states that
\begin{equation}
\label{eqn:hard}
R^- \in \CV^{-,T} \backslash \CS^{-,T}  \quad \Rightarrow \quad \exists \phi \in J^+ \cap \tUU^{+,T} \text{ s.t. } \Big \langle \phi, R^- \Big \rangle \neq 0 
\end{equation}
Given $R^- \in \CV^{-,T}$, the fact that $R^- \notin \CS^{-,T}$ implies via \eqref{eqn:dual 3} that there exists 
\begin{equation}
\label{eqn:choose phi}
\phi = e_{i_1,l_1} \dots e_{i_{n-k},l_{n-k}} W_{p,\bk,d-l_1-\dots-l_{n-k}} \in J^+
\end{equation}
for some $\b0 \leq \bk \leq \bn$ (we will henceforth denote $n = |\bn|$ and $k = |\bk|$), $p \in (\BK^\times)^n/\BK^\times$, $i_1,\dots,i_{n-k} \in I$ and $l_1,\dots,l_{n-k} \in \BZ$ such that 
\begin{equation}
\label{eqn:phi is 0}
\Big \langle \phi, R^- \Big \rangle \neq 0
\end{equation}
Moreover, we may assume the word
\begin{equation}
\label{eqn:the word}
v = \left[i_1^{(l_1)} \dots i_{n-k}^{(l_{n-k})} \right]
\end{equation}
to be non-increasing, and maximal such that \eqref{eqn:phi is 0} holds. We may also write
\begin{equation}
\label{eqn:choose w}
W_{p,\bk,d-l_1-\dots-l_{n-k}} =  \mathop{\sum_{i_{n-k+1}, \dots, i_n \in I}}_{l_{n-k+1}, \dots, l_n \in \BZ} \text{coefficient} \cdot e_{i_{n-k+1},l_{n-k+1}} \dots e_{i_n, l_n}
\end{equation}
where we only require in all terms of the sum above that
\begin{equation}
\label{eqn:close}
l_{n-k+1}, \dots, l_n \text{ are within a global constant away from their average}
\end{equation}
We can accomplish this because Remark \ref{rem:shift} for the shift vector $(1,\dots,1) \in \nn$ allows us to obtain all the elements $W_{p,\bk,\bullet}$ from finitely many values of $\bullet$. Then formula \eqref{eqn:phi is 0} implies that there are $i_1,\dots,i_n \in I$, $l_1\dots,l_n \in \BZ$ as above such that
$$
0 \neq \Big \langle  e_{i_1,l_1} \dots e_{i_n,l_n}, R^- \Big \rangle
$$
If we write
\begin{equation}
\label{eqn:r with coefficients}
R^- = \dots + \text{coefficient} \cdot z_1^{-k_1} \dots z_n^{-k_n} + \dots
\end{equation}
(with the variables relabeled as in \eqref{eqn:relabeling}), then we conclude that
\begin{align*}
0 &\neq \int_{|z_1| \gg \dots \gg |z_n|} \frac {\dots + \text{coefficient} \cdot z_1^{l_1-k_1} \dots z_n^{l_n-k_n}+ \dots}{\prod_{1\leq a < b \leq n} \zeta_{i_bi_a} \left( \frac {z_b}{z_a} \right)} \prod_{a=1}^n Dz_a \\
&= \int_{|z_1| \gg \dots \gg |z_n|} \frac {\dots + \text{coefficient} \cdot z_1^{l_1-d_1} \dots z_n^{l_n-d_n}+ \dots}{\prod_{1\leq a < b \leq n} \left( \alpha_{i_bi_a} + \sum_{s\geq 1} \text{coefficient}\cdot \frac {z_b^s}{z_a^s} \right)} \prod_{a=1}^n Dz_a
\end{align*}
where $d_1,\dots,d_n$ are related to $k_1,\dots,k_n$ via formula \eqref{eqn:d and k}. The maximality of the word \eqref{eqn:the word} implies that
\begin{equation}
\label{eqn:are equal}
(l_1,\dots,l_{n-k}) = (d_1,\dots,d_{n-k})
\end{equation}
for some monomial which appears with non-zero coefficient in $R^-$ (see \eqref{eqn:r with coefficients}).

\medskip

\noindent \textbf{Choice of $T$:} Let $m$ and $M$ be sufficiently large natural numbers, and define
\begin{equation}
\label{eqn:choice of t}
T = \left \{w = \left[i_1^{(d_1)} \dots i_n^{(d_n)} \right] \in \CW_{\text{non-inc}} \text{ such that }\deg w = (\bn, d) \text{ and } \right. 
\end{equation}
$$
\left. \sum_{a \in A} d_a \geq -M|A| + m|A|^2 + \sum_{A \ni a < b \notin A} (s_{i_ai_b} + s_{i_bi_a}) , \forall A \subseteq \{1,\dots,n\} \right\} 
$$
(since $M$ can be arbitrarily large, this would ensure the fact that any finite set of non-increasing words can be contained in $T$). The reason we add $\sum_{A \ni a < b \notin A} (s_{i_ai_b} + s_{i_bi_a})$ to the right-hand side of the inequality above is the straightforward fact (which we leave as an exercise to the interested reader) that the inequality on the second line of \eqref{eqn:choice of t} holds for the leading word \eqref{eqn:leading word} of a monomial \eqref{eqn:monomial} if and only if it holds for the analogous word associated to any other order of the variables in the monomial. Then \eqref{eqn:are equal} and the fact that $R^- \in \CV^{-,T}$ imply that
\begin{equation}
\label{eqn:ineq final}
\sum_{a \in A} l_a \geq -M|A| + m|A|^2 +  \sum_{A \ni a < b \notin A} (s_{i_ai_b} + s_{i_bi_a})
\end{equation}
where $A = B$ or $A = B \sqcup \{n-k+1,\dots,n\}$, for arbitrary $B \subseteq \{1,\dots,n-k\}$.

\medskip

\noindent We will use inequality \eqref{eqn:ineq final} to show that $\phi$ of \eqref{eqn:choose phi} lies in $\tUU^{+,T}$ (which would conclude the proof of \eqref{eqn:hard}, and with it Claim \ref{claim:claim}). Let $b \gg \max(\beta(1),\dots,\beta(n))$, with the $\beta$'s as in \eqref{eqn:straighten}. For any $i_1,\dots,i_n \in I$, $l_1,\dots,l_n \in \BZ$ which appear in \eqref{eqn:choose phi} and \eqref{eqn:choose w}, consider the largest index $x \in \{0,\dots,n-k\}$ such that
\begin{equation}
\label{eqn:final ineq}
l_x < l_{x+1} - b
\end{equation}
(we make the convention that $l_0 = -\infty$). We may use relation \eqref{eqn:straighten} to write 
$$
e_{i_{x+1},l_{x+1}} \dots e_{i_n,l_n} = \sum_{\left[j_{x+1}^{(l'_{x+1})} \dots j_n^{(l'_n)} \right] \in \CW_{\text{non-inc}}} \text{coefficient} \cdot e_{j_{x+1},l'_{x+1}} \dots e_{j_n,l'_n} 
$$
We need to make two observations about the $l'$'s that appear in the formula above.

\medskip

\begin{itemize}[leftmargin=*]

\item All the numbers $l'_{x+1},\dots,l'_n$ are within a global constant away from their average

\medskip

\item The large difference between $l_x$ and $l_{x+1}$ ensures that all concatenated words
\begin{equation}
\label{eqn:concatenated words}
\left[ i_1^{(l_1)} \dots i_x^{(l_x)} j_{x+1}^{(l'_{x+1})} \dots j_n^{(l'_n)} \right]
\end{equation}
which arise in the procedure above are non-increasing (recall that the word \eqref{eqn:the word} was non-increasing to begin with, and thus so are all of its prefixes)

\end{itemize}

\medskip

\noindent Thus, it remains to show that the concatenated words that appear in \eqref{eqn:concatenated words} are in $T$. By \eqref{eqn:ineq final}, we have 
\begin{equation}
\label{eqn:last 1}
\sum_{a \in B} l_a + \sum_{a \in C} l'_a \geq -M|A| + m|A|^2 +\sum_{A \ni a < b \notin A} (s_{i_ai_b} + s_{i_bi_a})
\end{equation}
where $A = B \sqcup C$ for arbitrary $B \subseteq \{1,\dots,x\}$, but $C = \varnothing$ or $C = \{x+1,\dots,n\}$. If the defining property of $T$ were to be violated, then we would have 
\begin{equation}
\label{eqn:last 2}
\sum_{a \in B} l_a + \sum_{a \in C} l'_a < -M|A| + m|A|^2 + \sum_{A \ni a < b \notin A} (s_{i_ai_b} + s_{i_bi_a})
\end{equation}
where $A = B \sqcup C$ for some $B \subseteq \{1,\dots,x\}$, but $C$ a \underline{proper} subset of $\{x+1,\dots,n\}$. We claim that \eqref{eqn:last 1} and \eqref{eqn:last 2} are incompatible (for $m$ chosen large enough compared to the global constant mentioned in the first bullet above, and the $s_{ij}$'s). Indeed, letting $\mu$ be the average from the first bullet above, relation \eqref{eqn:last 1} implies
\begin{align}
&\sum_{s \in B} l_s \geq -M |B|+m|B|^2 - \varepsilon \label{eqn:last 3} \\
&\sum_{s \in B} l_s + \mu(n-x) \geq -M(|B|+n-x) + m(|B|+n-x)^2 - \varepsilon \label{eqn:last 4}
\end{align}
for some sufficiently large global constant $\varepsilon$. Meanwhile, \eqref{eqn:last 2} implies
\begin{equation}
\label{eqn:last 5}
\sum_{s \in B} l_s + \mu  y < -M (|B|+y)+m(|B|+y)^2 + \varepsilon 
\end{equation}
where $y = |C|$ lies in $\{1,\dots,n-x-1\}$. Subtracting \eqref{eqn:last 3} from \eqref{eqn:last 5} yields
$$
\mu < -M+m(y +2|B|) + 2\varepsilon
$$
and subtracting \eqref{eqn:last 5} from \eqref{eqn:last 4} yields 
$$
\mu > -M + m(n-x+y + 2|B|) - 2 \varepsilon
$$
The two inequalities above are incompatible with each other if $m$ is chosen sufficiently larger than $\varepsilon$, thus yielding the desired contradiction.

\end{proof}

\end{proof}

\end{document}